\documentclass[12 pt]{elsarticle}

\usepackage{amssymb,amsmath,amscd,graphicx,
latexsym,amsthm,verbatim,color}
\usepackage[mathscr]{eucal}
\usepackage[all,cmtip]{xy}
\usepackage{pb-diagram}

\usepackage{multido}

\theoremstyle{plain}
\newtheorem{theorem}{Theorem}[section]
\newtheorem*{teorema}{Theorem}
\newtheorem{proposition}[theorem]{Proposition}
\newtheorem{corollary}[theorem]{Corollary}
\newtheorem{lemma}[theorem]{Lemma}

\numberwithin{equation}{section}

\theoremstyle{definition}
\newtheorem{application}[theorem]{Application}
\newtheorem{definition}[theorem]{Definition}

\newtheorem{remark}[theorem]{Remark}
\newtheorem{example}[theorem]{Example}

\newtheorem{notation}[theorem]{Notation}

\newcommand{\pe}{{{\bf{P}}}(\mathcal{E})}
\newcommand{\Ox}{\mathcal{O}_{X}}
\newcommand{\Ope}{\mathcal{O}_{{\bf{P}}(\mathcal{E})}}

\newcommand*{\longhookrightarrow}{\ensuremath{\lhook\joinrel\relbar\joinrel\rightarrow}} 
\newcommand{\sbl}{\vskip 3pt}

\begin{document}

\title{Okounkov bodies on projectivizations of rank two toric vector bundles}

\author{Jos\'e Luis Gonz\'alez\fnref{fn1}}
\ead{jgonza@umich.edu}

\fntext[fn1]{The author was partially supported by the NSF under grant DMS-0502170 and by the David and Lucile Packard Foundation under M. Musta\c{t}\u{a}'s fellowship.}
\address{Department of Mathematics, University of Michigan,
Ann Arbor, MI 48109, USA}




\setcounter{tocdepth}{1}


\begin{abstract}
The global Okounkov body of a projective variety is a closed convex cone that encodes asymptotic information about every big line bundle on the variety. In the case of a rank two toric vector bundle $\mathcal{E}$ on a smooth projective toric variety, we use its Klyachko filtrations to give an explicit description of the global Okounkov body of $\pe$. In particular, we show that this is a rational polyhedral cone and that $\pe$ is a Mori dream space. 
\end{abstract}

\maketitle

\vspace{0.5cm}

\tableofcontents

\vfill


\section{Introduction}

In his work on log-concavity of multiplicities, e.g. \cite{Okounkov96}, \cite{Okounkov03}, A. Okounkov introduced a procedure to associate convex bodies to linear systems on projective varieties.
This construction was systematically studied by R. Lazarsfeld and M. Musta\c{t}\u{a} in the case of big line bundles in \cite{CBALS}.
The construction of these \emph{Okounkov bodies} depends on a fixed flag of subvarieties and produces a convex compact set for each Cartier divisor on a projective variety.
The Okounkov body of a divisor encodes asymptotic invariants of the divisor's linear system, and it is determined solely by the divisor's  numerical equivalence class. Moreover, these bodies vary as fibers of a linear map defined on a closed convex cone as one moves in the space of numerical equivalence classes of divisors on the variety. As a consequence, one can expect to obtain results about line bundles by applying methods from convex geometry to the study of these Okounkov bodies.
 
Let us consider an $n$-dimensional projective variety $X$ over an algebraically closed field, endowed with a flag $X_{\bullet} \colon X = X_{n} \supseteq \cdots \supseteq X_{0}=\{pt\}$, where $X_{i}$ is an $i$-dimensional subvariety that is nonsingular at the point $X_{0}$. In \cite{CBALS}, Lazarsfeld and Musta\c{t}\u{a} established the following: \\
\textbf{(a)} 
For each big rational numerical divisor class $\xi$ on X, Okounkov's construction yields a convex compact set $\Delta(\xi)$ in ${\bf{R}}^n$, now called the \emph{Okounkov body} of $\xi$, whose Euclidean volume satisfies
\[ \textnormal{vol}_{{\bf{R}}^n}\big( \Delta(\xi) \big) \ = \ \frac{1}{n!} \cdot \textnormal{vol}_X(\xi).\]
The quantity $\textnormal{vol}_X(\xi)$ on the right is the \emph{volume} of the rational class $\xi$, which is defined by extending the definition of the volume of an integral Cartier divisor $D$ on $X$, namely, 
\[  
\textnormal{vol}_{X}(D) \ =_{\text{def}} \ \lim_{m \to \infty} \frac{ h^{0}(X,\Ox(mD))}{m^n / n!}.
\]
We recall that the volume is an interesting invariant of big divisors which plays an important role in several recent developments in higher dimensional geometry. For basic properties of volumes we refer to \cite{PAG}.  \\
\textbf{(b)} 
Moreover, there exists a closed convex cone $\Delta(X) \subseteq {\bf{R}}^{n} \times N^{1} (X)_{{\bf{R}}}$ characterized by the property that in the diagram 
	\[
\xymatrix{
\ \Delta(X) \ \ar@/_/[dr] \ar@{^{(}->}[rr] & & \ {\bf{R}}^{n} \times N^{1} (X)_{{\bf{R}}} \ar@/^/[ld] \ \\
& N^1(X)_{{\bf{R}}},  }
	\]
the fiber $\Delta(X)_\xi \subseteq {\bf{R}}^n \times\{ \xi \} = {\bf{R}}^n$ of $\Delta(X)$ over any big class $\xi \in N^1(X)_{{\bf{Q}}}$ is $\Delta(\xi)$. This is illustrated schematically in Figure 1. $\Delta(X)$ is called the \emph{global Okounkov body} of $X$.



\begin{center}
\includegraphics[width=3in]{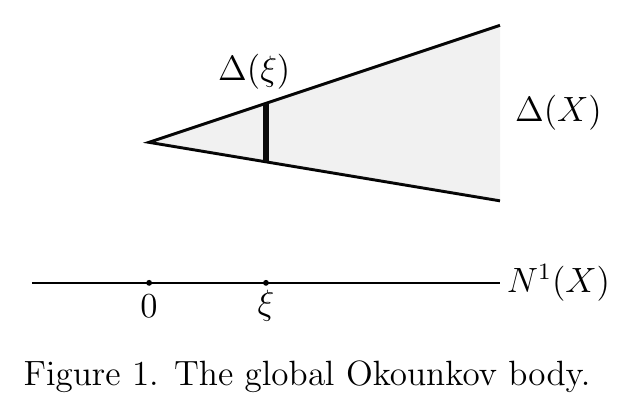}
\end{center}


Lazarsfeld and Musta\c{t}\u{a} have used this theory to reprove and generalize results about volumes of divisors, including Fujita's Approximation Theorem. Using $\bf{(b)}$, they can additionally give alternative proofs of properties of the volume function $\textnormal{vol}_{X} \colon \operatorname{Big}(X) \rightarrow {\bf{R}}$, defined in the set of big classes of ${\bf{R}}$-divisors. For example, it follows that $\textnormal{vol}_{X}$ is of class $\mathcal{C}^{1}$ and satisfies the log-concavity relation  
\[ 
\textnormal{vol}_{X}(\xi + \xi^\prime)^{1/n}  \ \ge \ \textnormal{vol}_{X}(\xi)^{1/n} + \textnormal{vol}_{X}(\xi^\prime)^{1/n},
\]
for any two big classes $\xi, \xi^\prime \in N^1(X)_{{\bf{R}}}$.
It is worth mentioning that in \cite{KK}, K. Kaveh and A. Khovanskii use a similar procedure to associate convex bodies to finite dimensional subspaces of the function field $\operatorname{K}(X)$ of a variety $X$. In their work, they use the bridge between convex and algebraic geometry provided by their construction to obtain results in both areas.

The explicit description of Okounkov bodies in concrete examples can be rather difficult. One easy case is that of smooth projective toric varieties. If $D$ is an invariant divisor on such a variety, and if the flag consists of invariant subvarieties, then it is shown in \cite{CBALS} that the Okounkov body of $D$ is the polytope $P_{D}$ that one usually associates to $D$ in toric geometry, up to a suitable translation. 

In this paper, we are interested in describing the Okounkov bodies of the divisors on the projectivization $\pe$ of a toric vector bundle $\mathcal{E}$ on the smooth projective toric variety $X$. Such vector bundles were described by A. Klyachko in \cite{Klyachko} in terms of certain filtrations of a suitable vector space, and they have been the focus of some recent activity, e.g. \cite{PTVB}, \cite{Payne}, \cite{Payne2}. As we will see, these filtrations can be used to compute the sections of all line bundles on $\pe$.
For our main result, we restrict to the case of rank two toric vector bundles, where the Klyachko filtrations are considerably simpler. Using the data from the filtrations, we construct a flag of torus invariant subvarieties on $\pe$ and produce finitely many linear inequalities defining the global Okounkov body of $\pe$ with respect to this flag. In particular, we see that this is a rational polyhedral cone. As an application we show that the Cox ring of $\pe$ is finitely generated, giving an alternative proof in a special case of a result of Knop (see \cite{Knop} and \cite{Suess}). 

This article is organized as follows. In \S \ref{section.notation} we review Okounkov's construction and Klyachko's description of toric vector bundles. Next, in \S \ref{section.filtrations} we describe the Klyachko filtrations for tensor products and Schur functors of toric vector bundles of arbitrary rank. This allows us to compute the global sections of any line bundle on $\pe$, a fact of independent interest.
The main point in the description of Okounkov bodies is the computation of the vanishing orders of all sections of certain line bundles along the subvarieties in a fixed flag. In \S \ref{section.vanishing.orders} we construct a flag of invariant subvarieties $Y_{\bullet}$ in $\pe$, and show that in this setting it is enough to compute the vanishing orders of a special collection of isotypical sections with respect to the torus action. Finally, in \S \ref{section.okounkov.body} we describe the global Okounkov body of $\pe$ with respect to the flag $Y_{\bullet}$ and prove that the variety $\pe$ is a Mori dream space. In the last section, we give some examples to illustrate our main result.

\subsection*{Acknowledgments}
I would like to thank W. Fulton, M. Hering, R. Lazarsfeld and S. Payne, from whom I learned some of the topics relevant to this article. I express my gratitude to the anonymous referee for his/her careful reading of this paper and his/her useful suggestions. I would especially like to thank M. Musta\c{t}\u{a} for many insightful conversations and for his comments on earlier versions of this manuscript.





\section{Okounkov bodies and toric vector bundles} \label{section.notation}

In this section, we review the construction of Okounkov bodies and Klyachko's classification of toric vector bundles. Unless explicitly stated otherwise, the notation introduced here will be used throughout the paper. All our varieties are assumed to be defined over a fixed algebraically closed field $k$. By a \emph{divisor} on a variety $Z$ we always mean a Cartier divisor on $Z$. We denote the group of numerical equivalence classes of divisors on $Z$ by $N^{1}(Z)$, and we denote the spaces $N^{1}(Z) \otimes {\bf{Q}}$ and $N^{1}(Z) \otimes {\bf{R}}$ by $N^{1}(Z)_{{\bf{Q}}}$ and $N^{1}(Z)_{{\bf{R}}}$, respectively.

By a \emph{line bundle} on a variety $Z$, we mean a locally free sheaf of rank one on $Z$. We follow the convention that the \emph{geometric vector bundle} associated to the locally free sheaf $\mathcal{F}$ is the variety ${\bf{V}}(\mathcal{F}) =  \text{\textbf{Spec}} \, \bigoplus_{m \geq 0} Sym^{m}\mathcal{F}^{\vee}$, whose sheaf of sections is $\mathcal{F}$. Also, by the fiber of $\mathcal{F}$ over a point $z \in Z$, we mean the fiber over $z$ of the projection $f \colon {\bf{V}}(\mathcal{F}) \rightarrow Z$. Lastly, by the \emph{projectivization} ${\bf{P}}(\mathcal{F})$ of $\mathcal{F}$, we mean the projective bundle $\text{\textbf{Proj}} \, \bigoplus_{m \geq 0} Sym^{m}\mathcal{F}$ over $Z$. This bundle is endowed with a projection $\pi \colon {\bf{P}}(\mathcal{F}) \rightarrow Z $ and an invertible sheaf $\mathcal{O}_{{\bf{P}}(\mathcal{F})}(1)$
(see II.7 in \cite{Hartshorne}).


\subsection{Okounkov bodies} \label{okounkov.procedure}
Let us consider a normal $l$-dimensional variety $Z$ with a fixed flag $Y_{\bullet} \colon Z=Y_{l} \supseteq \cdots \supseteq Y_{0}$, where each $Y_{i}$ is an $i$-dimensional normal subvariety that is nonsingular at the point $Y_{0}$. Given a big divisor $D$ on $Z$, we will describe a procedure to assign a compact convex set with nonempty interior $\Delta_{Y_{\bullet}}(D)$ in ${\bf{R}}^{l}$ to $D$.
First, given any divisor $F$ on $Z$ and a nonzero section $s=s_{l} \in H^{0}(Y_{l},\mathcal{O}_{Y_{l}}(F))$, we can associate an $l$-tuple of nonnegative integers $\nu_{Y_{\bullet},F}(s) = (\nu_{1}(s),\ldots,\nu_{l}(s))$ to $s$ as follows. By restricting to a neighborhood of $Y_{0}$ we can assume that each $Y_{i}$ is smooth. We define $\nu_{1}(s)$ to be the vanishing order $\operatorname{ord}_{Y_{l-1}}(s)$ of $s$ along $Y_{l-1}$. Then, $s$ determines a section $\widetilde{s_{l}} \in H^{0}(Y_{l},\mathcal{O}_{Y_{l}}(F) \otimes \mathcal{O}_{Y_{l}}(-\nu_{1}(s)Y_{l-1}))$ that does not vanish along $Y_{l-1}$. By restricting, we get a nonzero section $s_{l-1} \in H^{0}(Y_{l-1},\mathcal{O}_{Y_{l}}(F)|_{Y_{l-1}} \otimes \mathcal{O}_{Y_{l}}(-\nu_{1}(s)Y_{l-1})|_{Y_{l-1}})$, and we iterate this procedure. More precisely, assume that we have defined nonnegative integers $\nu_{1}(s),\ldots,\nu_{h}(s)$, and nonzero sections $s_{l} \in H^{0}(Y_{l},\mathcal{O}_{Y_{l}}(F)),\ldots,s_{l-h} \in H^{0} ( Y_{l-h},\mathcal{O}_{Y_{l}}(F)|_{Y_{l-h}} \otimes \bigotimes_{i=1}^{h} \mathcal{O}_{Y_{l-i+1}}(-\nu_{i}(s)Y_{l-i})|_{Y_{l-h}})$, for some nonnegative integer $h < l$. We define $\nu_{h+1}(s)$ as the vanishing order $\operatorname{ord}_{Y_{l-h-1}}(s_{l-h})$ of $s_{l-h}$ along $Y_{l-h-1}$; then, $s_{l-h}$ determines a section 
\[
\widetilde{s_{l-h}} \in H^{0} ( Y_{l-h},\mathcal{O}_{Y_{l}}(F)|_{Y_{l-h}} \otimes \bigotimes_{i=1}^{h+1} \mathcal{O}_{Y_{l-i+1}}(-\nu_{i}(s)Y_{l-i})|_{Y_{l-h}})
\]
that does not vanish along $Y_{l-h-1}$; and finally, by restricting, we get a nonzero section 
\[
s_{l-h-1} \in H^{0} ( Y_{l-h-1},\mathcal{O}_{Y_{l}}(F)|_{Y_{l-h-1}} \otimes \bigotimes_{i=1}^{h+1} \mathcal{O}_{Y_{l-i+1}}(-\nu_{i}(s)Y_{l-i})|_{Y_{l-h-1}}).
\]
We repeat this procedure until we obtain nonnegative integers $\nu_{1}(s),\ldots,\nu_{l}(s)$. This construction gives us a function
\begin{align*}
\begin{matrix}
\nu_{Y_{\bullet}} = \nu_{Y_{\bullet},F} \colon&H^{0}(Z,\mathcal{O}_{Z}(F)) \smallsetminus \{ 0 \}&\longrightarrow&{\bf{Z}}^{l}\\
&s&\longmapsto&(\nu_{1}(s),\ldots,\nu_{l}(s)).
\end{matrix}
\end{align*}
We denote the image of $\nu_{Y_{\bullet},F}$ by either $\nu(F)$ or $\nu(\mathcal{O}_{Z}(F))$. The function $\nu_{Y_{\bullet}}$ satisfies the following valuation-like properties:
\begin{itemize}
\item For any nonzero sections $s_{1},s_{2} \in H^{0}(Z,\mathcal{O}_{Z}(F))$, we have that $\nu_{Y_{\bullet},F}(s_{1}+s_{2}) \geq_{lex} \linebreak {\operatorname{min}}_{\geq_
{lex}}  \{\nu_{Y_{\bullet},F}(s_{1}),\nu_{Y_{\bullet},F}(s_{2}) \}$, where $\geq_{lex}$ denotes the lexicographic order in ${\bf{Z}}^{l}$.
\item For any divisors $F_{1}$ and $F_{2}$ in $Z$, and nonzero sections $s_{1} \in H^{0}(Z,\mathcal{O}_{Z}(F_{1}))$ and $s_{2} \in H^{0}(Z,\mathcal{O}_{Z}(F_{2}))$, we have $\nu_{Y_{\bullet},F_{1}+F_{2}}(s_{1} \otimes s_{2}) = \nu_{Y_{\bullet},F_{1}}(s_{1}) + \nu_{Y_{\bullet},F_{2}}(s_{2})$.
\end{itemize}


\begin{remark} \label{Lemma:LazMus}
If $W$ is a finite dimensional subspace of $H^{0}(Z,\mathcal{O}_{Z}(F))$, then the number of valuation vectors arising from nonzero sections in $W$ is equal to the dimension of $W$. For example, when $Z$ is complete, $\nu(F)$ is a finite set with cardinality $\operatorname{dim}_{k} H^{0}(Z,\mathcal{O}_{Z}(F))$. A more general statement is proven by Lazarsfeld and Musta\c{t}\u{a} as Lemma 1.3 in \cite{CBALS}.
\end{remark}

Finally, $\Delta_{Y_{\bullet}}(F) = \Delta_{Y_{\bullet}}(\mathcal{O}_{Z}(F))$ is defined to be the following closed convex hull in ${\bf{R}}^{l}$: 
\[
\Delta_{Y_{\bullet}}(F) = \overline{ \operatorname{Conv} \Big{(} \bigcup_{m \in {\bf{Z}}^{+}} \frac{1}{m} \nu ( m F ) \Big{)} }.
\]
We will denote the set $\Delta_{Y_{\bullet}}(F)$ simply by either $\Delta(F)$ or $\Delta(\mathcal{O}_{Z}(F))$ whenever the corresponding flag is understood. In \cite{CBALS}, Lazarsfeld and Musta\c{t}\u{a} proved that when $Z$ is a projective variety and $D$ is a big divisor, $\Delta(D)$ is a compact convex subset of ${\bf{R}}^{l}$ with nonempty interior, i.e. a \emph{convex body}. In this case, $\Delta(D)$ is called the \emph{Okounkov body} of D. 

Since $\Delta_{Y_{\bullet}}(mF) = m \Delta_{Y_{\bullet}}(F)$ for any divisor $F$ and any $m \in {\bf{Z}}^{+}$, this definition extends in a natural way to associate an Okounkov body to any big divisor with rational coefficients. As it turns out, the outcome depends only on the numerical equivalence class of the divisor. We refer to \textbf{(a)} and \textbf{(b)} in the introduction for some of the main properties of this construction, including the existence of the \emph{global Okounkov body} of a projective variety $Z$. This global Okounkov body is a closed convex cone $\Delta(Z) \subseteq {\bf{R}}^{l} \times N^{1} (Z)_{{\bf{R}}}$ characterized by the property that the fiber of the second projection over any big class $D \in N^{1} (Z)_{{\bf{Q}}}$ is the Okounkov body $\Delta(D)$. For proofs of these assertions, as well as of \textbf{(a)} and \textbf{(b)}, we refer to \cite{CBALS}. 

\begin{example}
Let $Z = {\bf{P}}^l$ with homogeneous coordinates $z_{0},\ldots,z_{l}$. Let $Y_{\bullet}$ be the flag of linear subspaces defined by $Y_{i} = \{ z_{1} = \cdots = z_{l-i} = 0\}$ for each $i$. If $|D|$ is the linear system of hypersurfaces of degree $m$, then $\nu_{Y_{\bullet},D}$ is the lexicographic valuation determined on monomials of degree $m$ by 
\[
\nu_{Y_{\bullet}}(z_{0}^{\alpha_{0}} \cdots z_{l}^{\alpha_{l}}) =  (\alpha_{1},\ldots,\alpha_{l}) , 
\]
and the Okounkov body $\Delta(D)$ is the simplex 
\[
\Delta(D) = \{(\lambda_{1},\ldots,\lambda_{l}) \in {\bf{R}}^{l} \ | \ \lambda_{1} \geq 0, \ldots, \lambda_{l} \geq 0, \sum_{i=1}^{l} \lambda_{i} \leq m \}.
\] 
\end{example}


\subsection{Toric vector bundles and Klyachko filtrations} \label{subsection.toric}
Let $X$ be an $n$-dimensional toric variety corresponding to a fan $\Delta$ in the lattice $N$. We denote the algebraic torus acting on $X$ by $T$, and the character lattice $\operatorname{Hom}(N,{\bf{Z}})$ of $T$ by $M$. Thus, $T=\operatorname{Spec}k[M] = \operatorname{Spec}k[ \, \chi^{u} \, | \, u \in M \, ]$ and $X$ has an open covering given by the affine toric varieties $U_{\sigma} = \operatorname{Spec}k[ \, \sigma^{\vee} \cap M \, ]$ corresponding to the cones $\sigma \in \Delta$. We denote the rays in $\Delta$ by $\rho_{1},\ldots,\rho_{d}$. For each ray $\rho_{j}$, we denote its primitive lattice generator by $v_{j}$ and its associated codimension one torus invariant subvariety by $D_{j}$. Let $t_{0}$ denote the unit element of the torus. For a detailed treatment of toric varieties we refer to \cite{Fulton}. 

If $T$ acts on a vector space $V$ in such a way that each element of $V$ belongs to a finite dimensional $T$-invariant subspace, we get a decomposition $V = \bigoplus_{u \in M} V_{u}$, where $V_{u} =_{def} \{ v \in V  \ | \  tv = \chi^{u} (t) v \textnormal{ for each $t \in T$}  \}$. The spaces $V_{u}$ and their elements are called isotypical summands and isotypical elements, respectively. This motivates the use of the following terminology. If $T$ acts on the space of sections of a vector bundle on some variety, we say that a section $s$ is $T$-\emph{isotypical} if there exists $u \in M$ such that $ts=\chi^{u}(t)s$ for each $t \in T$. Likewise, if $T$ acts on a variety, we say that a rational function $f$ on the variety is $T$-\emph{isotypical} if there exists $u \in M$ such that $tf=\chi^{u}(t)f$ for each $t \in T$, i.e. the domain $\operatorname{dom}(f)$ of $f$ is $T$-invariant and $(tf)(z) =_{def} f(t^{-1}z)  = \chi^{u}(t)f(z)$ for each $z \in \operatorname{dom}(f)$ and each $t \in T$. When $T$ acts algebraically on an affine variety $Z$, the induced action of $T$ on $H^{0}(Z, \mathcal{O}_{Z})$ satisfies the above finiteness condition, and we get a decomposition $H^{0}(Z, \mathcal{O}_{Z}) = \bigoplus_{u \in M} H^{0}(Z, \mathcal{O}_{Z})_{u}$ as before (see I.6.3 in \cite{LePotier}). 

A \emph{toric vector bundle} on the toric variety $X$ is a locally free sheaf $\mathcal{E}$ together with an action of the torus $T$ on the variety ${\bf{V}}(\mathcal{E})$, such that the projection $f \colon {\bf{V}}(\mathcal{E}) \rightarrow X$  is equivariant and $T$ acts linearly on the fibers of $f$. In general, if $\mathcal{E}$ is a toric vector bundle, ${\bf{V}}(\mathcal{E})$ and $\pe$ need not be toric varieties. Given any $T$-invariant open subset $U$ of $X$, there is an induced action of $T$ on $H^{0}(U,\mathcal{E})$, defined by the equation 
\[
(t \cdot s)(x)=_{\operatorname{def}} t (s( t^{-1} x )),
\]
for any $s \in H^{0}(U,\mathcal{E})$, $t \in T$ and $x \in X$. This action induces a direct sum decomposition 
\[
H^{0}(U,\mathcal{E}) = \bigoplus_{u \in M} H^{0}(U,\mathcal{E})_{u},
\]
where $H^{0}(U,\mathcal{E})_{u}= \{s \in H^{0}(U,\mathcal{E}) \ | \ t \cdot s = \chi^{u} (t) s \textnormal{ for each $t \in T$}  \}$, as before.

\begin{example}
For each torus invariant Cartier divisor $D$ on a toric variety $X$, the line bundle $\Ox(D)$ has a natural toric vector bundle structure. For each $w \in M$, the isotypical summands in the decomposition of $H^{0}(U_{\sigma},\mathcal{O}_{X}(\operatorname{div}\chi^{w}))$ over $U_{\sigma}$ are given by
\begin{equation*} 
H^{0}(U_{\sigma},\mathcal{O}_{X}(\operatorname{div}\chi^{w}))_{u}=
	\begin{cases}
		k \chi^{-u} & \text{if $w-u \in \sigma^{\vee} \cap M$},\\
		0& \text{otherwise},
	\end{cases}
\end{equation*}
for each cone $\sigma \in \Delta$ and each $u \in M$.
\end{example}

A toric vector bundle over an affine toric variety is equivariantly isomorphic to a direct sum of toric line bundles (see Proposition 2.2 in \cite{Payne2}). Every line bundle $\mathcal{L}$ on $X$ admits a T-equivariant structure, and choosing one such structure is equivalent to choosing a torus invariant divisor $D$ such that $\mathcal{L} \cong \Ox(D)$. The classification of toric vector bundles of higher rank is considerably more complicated.

Let $E$ be the fiber over $t_{0}$ of the toric vector bundle $\mathcal{E}$ on $X$. For each ray $\rho_{j} \in \Delta$ and each $u \in M$, the evaluation map at $t_{0}$ gives an inclusion $H^{0}(U_{\rho_{j}},\mathcal{E})_{u} \longhookrightarrow E$. If $u, u' \in M$ satisfy $\langle u ,v_{j} \rangle \geq \langle u' ,v_{j} \rangle$, then 
\[
\operatorname{Im}(H^{0}(U_{\rho_{j}},\mathcal{E})_{u} \longhookrightarrow E) \subseteq \operatorname{Im} (H^{0}(U_{\rho_{j}},\mathcal{E})_{u'} \longhookrightarrow E).
\] 
Therefore the images of these maps depend only on $\langle u ,v_{j} \rangle$, or equivalently, only on the class of $u$ in $M/\rho_{j}^{\perp}\cap M \cong {\bf{Z}}$. Hence, we may denote the image of the map $H^{0}(U_{\rho_{j}},\mathcal{E})_{u} \longhookrightarrow E$ simply by $\mathcal{E}^{\rho_{j}}(\langle u ,v_{j} \rangle)$. Note that for each $u \in M$ the image of the evaluation map $H^{0}(X,\mathcal{E})_{u} \longhookrightarrow E$ is equal to $\mathcal{E}^{\rho_{1}}(\langle u ,v_{1} \rangle)\cap\cdots\cap\mathcal{E}^{\rho_{d}}(\langle u ,v_{d} \rangle) \subseteq E$.
The ordered collection of finite dimensional vector subspaces $\mathcal{E}^{\rho_{j}}=_{def}\{\mathcal{E}^{\rho_{j}}(i) \ | \ i \in {\bf{Z}} \}$ gives a decreasing filtration of $E$. The filtrations $\{ \mathcal{E}^{\rho_{j}} \ | \ j=1,\ldots,d \ \}$ are called the \emph{Klyachko filtrations} associated to $\mathcal{E}$. For each $\sigma \in \Delta$, by equivariantly trivializing $\mathcal{E}$ over the affine open subset $U_{\sigma}$ of $X$, one can show that there exists a decomposition $E= \bigoplus_{\bar{u} \in M/\sigma^{\perp} \cap M}E_{\bar{u}}$, such that $\mathcal{E}^{\rho}(i)=\sum_{\langle \bar{u},v_{\rho} \rangle \geq i}{E_{\bar{u}}}$, for each ray $\rho \subseteq \sigma$ and each $i \in {\bf{Z}}$. Klyachko proved in \cite{Klyachko} that the vector space $E$ together with these filtrations, satisfying the above compatibility condition, completely describe $\mathcal{E}$. More precisely, 
\begin{teorema}[Klyachko] \label{theorem.klyachko.filtrations} The category of toric vector bundles on the toric variety $X$ is equivalent to the category of finite dimensional $k$-vector spaces $E$ with collections of decreasing filtrations $\{\mathcal{E}^{\rho}(i) \ | \ i \in {\bf{Z}} \}$, indexed by the rays of $\Delta$, satisfying the following compatibility condition: For each cone $\sigma \in \Delta$, there is a decomposition $E= \bigoplus_{\bar{u} \in M/\sigma^{\perp} \cap M}E_{\bar{u}}$ such that
\[
\mathcal{E}^{\rho}(i)=\sum_{\langle \bar{u},v_{\rho} \rangle \geq i}{E_{\bar{u}}},
\]
for every ray $\rho \subseteq \sigma$ and every $i \in {\bf{Z}}$.
\end{teorema}

\begin{example} \label{example.klyachko}
Let $m_{1}D_{1}+\cdots+m_{d}D_{d}$ be a torus invariant Cartier divisor on $X$, for some $m_{1},$ $\ldots,$ $m_{d}\in {\bf{Z}}$. Denote by $\mathcal{L}$ the line bundle $\Ox(m_{1}D_{1}+\cdots+m_{d}D_{d})$ on $X$. The Klyachko filtrations for $\mathcal{L}$ are given by
\begin{equation*} 
\mathcal{L}^{\rho_{j}}(i)=
	\begin{cases}
		k& \text{if $i \leq m_{j}$},\\
		0& \text{if $i > m_{j}$},
	\end{cases}
\end{equation*}
for each ray $\rho_{j} \in \Delta$.
\end{example}

\begin{example} \label{tangent.bundle.definition}
The projective plane ${\bf{P}}^{2}$ can be represented as the toric variety associated to the fan in $N \otimes {\bf{R}}={\bf{R}}^{2}$ with maximal cones $\sigma_{1} = \langle v_{2},v_{3} \rangle$, $\sigma_{2} = \langle v_{3},v_{1} \rangle$ and $\sigma_{3} = \langle v_{1},v_{2} \rangle$, where $v_{1}=(1,0)$, $v_{2}=(0,1)$ and $v_{3}=(-1,-1)$. 
The tangent bundle $T_{{\bf{P}}^{2}}$ of ${\bf{P}}^{2}$ is naturally a toric vector bundle on ${\bf{P}}^{2}$. This bundle can be equivariantly trivialized as
\begin{gather*}
T_{{\bf{P}}^{2}}|_{U_{\sigma_{1}}} = \mathcal{O}_{{\bf{P}}^{2}}(D_{2})|_{U_{\sigma_{1}}} \oplus  \mathcal{O}_{{\bf{P}}^{2}}(D_{3})|_{U_{\sigma_{1}}}, \qquad
T_{{\bf{P}}^{2}}|_{U_{\sigma_{2}}} = \mathcal{O}_{{\bf{P}}^{2}}(D_{3})|_{U_{\sigma_{2}}} \oplus \mathcal{O}_{{\bf{P}}^{2}}(D_{1})|_{U_{\sigma_{2}}}, \\
T_{{\bf{P}}^{2}}|_{U_{\sigma_{3}}} = \mathcal{O}_{{\bf{P}}^{2}}(D_{1})|_{U_{\sigma_{3}}} \oplus \mathcal{O}_{{\bf{P}}^{2}}(D_{2})|_{U_{\sigma_{3}}}. 
\end{gather*}
It follows that the Klyachko filtrations associated to $T_{{\bf{P}}^{2}}$ are
\begin{align*}
\begin{matrix}
T_{{\bf{P}}^{2}}^{\rho_{1}}(i)=
	\begin{cases}
		E& \text{if $i \leq 0$},\\
		V_{1}& \text{if $i = 1$}, \\
		0& \text{if $i \geq 2$},
	\end{cases} & &
T_{{\bf{P}}^{2}}^{\rho_{2}}(i)=
	\begin{cases}
		E& \text{if $i \leq 0$},\\
		V_{2}& \text{if $i = 1$},\\
		0& \text{if $i \geq 2$},
	\end{cases} & &
T_{{\bf{P}}^{2}}^{\rho_{3}}(i)=
	\begin{cases}
		E& \text{if $i \leq 0$},\\
		V_{3}& \text{if $i = 1$},\\
		0& \text{if $i \geq 2$},
	\end{cases} 
\end{matrix}
\end{align*}
where $V_{1}$, $V_{2}$ and $V_{3}$ are distinct one-dimensional subspaces of the fiber $E$ of $T_{{\bf{P}}^{2}}$ over $t_{0}$.
\end{example}

For any toric vector bundle $\mathcal{E}$ over $X$ of rank at least two, we have an isomorphism $N^{1}(X) \oplus {\bf{Z}} = \operatorname{Pic} X \oplus {\bf{Z}} \cong \operatorname{Pic}\pe = N^{1}(\pe)$, which is induced by $(\mathcal{L},m)\mapsto \Ope (m) \otimes \pi^{*}\mathcal{L}$, where $\pi$ is the projection map $\pi \colon \pe \rightarrow X$.



\section{The Klyachko filtrations for tensor products and Schur functors} \label{section.filtrations}

As we reviewed in \S \ref{theorem.klyachko.filtrations}, Klyachko proved that the category of toric vector bundles on a toric variety $X$ is equivalent to the category of finite dimensional vector spaces endowed with a collection of filtrations that satisfy a certain compatibility condition. Klyachko's result allows us to carry out some explicit computations in this category, including the description of the space of sections of a toric vector bundle over any invariant open subset of $X$. 
Throughout this section $X$ denotes an arbitrary toric variety.
 
Using the notation introduced in \S \ref{section.notation}, each line bundle on $\pe$ is isomorphic to a line bundle of the form $\Ope(m) \otimes \pi^{*}(\mathcal{O}_{X}(D))$ for some $T$-invariant Cartier divisor $D$ on $X$. For such an isomorphic representative we have a toric vector bundle structure on $\pi_{*}\bigl(\Ope (m) \otimes \pi^{*} \mathcal{O}_{X}(D)\bigr)=(Sym^{m}\mathcal{E}) \otimes  \mathcal{O}_{X}(D)$, and we have
	\[
H^{0}(\pe,\Ope(m) \otimes \pi^{*}(\mathcal{O}_{X}(D))=H^{0}(X, ( Sym^{m}\mathcal{E} ) \otimes \mathcal{O}_{X}(D)).
\]
From this, we see that the Klyachko filtrations associated to tensor products and symmetric powers of toric vector bundles can be used to describe the spaces of global sections of line bundles on projectivized toric vector bundles.
The goal of this section is to provide appropriate descriptions of these filtrations for toric vector bundles of arbitrary rank. We present the filtrations for tensor products in Lemma \ref{Lemma:TensorProductTVB} and its Corollary \ref{tensor.product.several}. The filtrations for symmetric powers are given in Corollary \ref{corollary.klyachko.symmetric} to Lemma \ref{lemma.klyachko.schur}. More generally, in that lemma we describe the filtrations for any Schur functor, e.g. symmetric and wedge products. 

We start by presenting the filtrations for tensor products.

\begin{lemma}\label{Lemma:TensorProductTVB}
Let $\mathcal{E}$ and $\mathcal{F}$ be toric vector bundles on the toric variety $X$. Then the Klyachko filtrations for their tensor product $\mathcal{E} \otimes \mathcal{F}$ are given by
\[
\big{(} \mathcal{E} \otimes \mathcal{F} \big{)} ^{\rho}(i)=\sum_{i_{1}+i_{2}=i} \mathcal{E}^{\rho}(i_{1}) \otimes  \mathcal{F}^{\rho}(i_{2}),
\]
for each ray $\rho \in \Delta$ and each $i\in{\bf{Z}}$.
\end{lemma}
\begin{proof}
Since the filtration corresponding to the ray $\rho$ only depends on $U_{\rho}$, it suffices to consider the case when $X=U_{\rho}$ for some ray $\rho \in \Delta$. Hence we can assume that $\mathcal{E}$ and $\mathcal{F}$ equivariantly trivialize as
	\begin{align*}
		\mathcal{E}&=\Ox(d_{1}D_{\rho}) \oplus \Ox(d_{2}D_{\rho}) \oplus \cdots \oplus \Ox(d_{r}D_{\rho}) \\
		\mathcal{F}&=\Ox(e_{1}D_{\rho}) \oplus \Ox(e_{2}D_{\rho}) \oplus \cdots \oplus \Ox(e_{s}D_{\rho})
	\end{align*}
for some $d_{1},\ldots,d_{r},e_{1},\ldots,e_{s} \in {\bf{Z}}$. Now we note that
\[
\big( \Ox(d_{j_{1}}D_{\rho}) \otimes \Ox(e_{j_{2}}D_{\rho}) \big)^{\rho}(i)=\sum_{i_{1}+i_{2}=i}{\Ox(d_{j_{1}}D_{\rho})^{\rho}(i_{1}) \otimes \Ox(e_{j_{2}}D_{\rho})^{\rho}(i_2)} 
\]
for each $i \in {\bf{Z}}$ and each $j_{1} \in \{ 1,\ldots,r \}$ and $j_{2} \in \{ 1,\ldots,s \}$. Since the Klyachko filtrations for a direct sum are the direct sums of the filtrations for the summands, the result now follows.
\end{proof}

\begin{corollary} \label{tensor.product.several}
Let $\mathcal{E}_{1},\ldots,\mathcal{E}_{s}$ be toric vector bundles on the toric variety $X$. Then the Klyachko filtrations	for their tensor product $\mathcal{E}_{1} \otimes \cdots \otimes \mathcal{E}_{s}$ are given by
\[
\big{(} \mathcal{E}_{1} \otimes \cdots \otimes \mathcal{E}_{s} \big{)}^{\rho}(i)=\sum_{i_{1}+\cdots+i_{s}=i}\mathcal{E}_{1}^{\rho}(i_{1}) \otimes \cdots \otimes  \mathcal{E}_{s}^{\rho}(i_{s}),
\]
for each $i\in{\bf{Z}}$ and each ray $\rho \in \Delta$.
\end{corollary}
\begin{proof}
The conclusion follows from the previous lemma by induction on $s$.
\end{proof}

\begin{example} \label{tensor.product}
Let $\mathcal{E}$ be a toric vector bundle on the toric variety $X$, and let $D = m_{1}D_{1}+\cdots+m_{d}D_{d}$ be a torus invariant Cartier divisor on $X$, for some $m_{1},\ldots,m_{d}\in {\bf{Z}}$. Let us denote the fiber over $t_{0}$ of the line bundle $\Ox(D)$ by $G$. From the previous lemma and Example \ref{example.klyachko}, it follows that the Klyachko filtrations of $\mathcal{E} \otimes \Ox(D)$ are given by
\[
(\mathcal{E} \otimes \Ox(D))^{\rho_{j}}(i)= \mathcal{E}^{\rho_{j}}(i-m_{j}) \otimes G,
\]
for each $i \in {\bf{Z}}$ and each ray $\rho_{j} \in \Delta$.
\end{example}

In the following lemma we describe the Klyachko filtrations for the toric vector bundle obtained by applying a Schur functor to another toric vector bundle. As a corollary, we state the case of symmetric products, which will be used in \S \ref{section.okounkov.body}. For the definition and basic properties of Schur functors we refer to \S 6 in \cite{Fulton.Harris}.


\begin{lemma} \label{lemma.klyachko.schur}
Let $\mathcal{E}$ be a toric vector bundle on the toric variety $X$, and let $S_{\lambda}$ be the Schur functor associated to a Young tableau $\lambda$ with $m$ entries. Then the Klyachko filtrations for $S_{\lambda}\mathcal{E}$ are given by
\[
\big( S_{\lambda}\mathcal{E} \big)^{\rho}(i)=\sum_{i_{1}+ \cdots + i_{m}=i} \operatorname{Im} \big( \mathcal{E}^{\rho}(i_{1}) \otimes \cdots \otimes \mathcal{E}^{\rho}(i_{m}) \longrightarrow S_{\lambda}(E) \big),
\]
for each ray $\rho \in \Delta$ and each $i\in{\bf{Z}}$.
\end{lemma}
\begin{proof}
Since $S_{\lambda}\mathcal{E}$ is a quotient of $\mathcal{E}^{\otimes m}$, it follows that $( S_{\lambda}\mathcal{E} )^{\rho}(i)$ is the image of $(\mathcal{E}^{\otimes m})^{\rho}(i)$ under the natural map $E^{\otimes m} \rightarrow S_{\lambda}(E)$, for each ray $\rho$ and each $i \in {\bf{Z}}$. Now, the result follows at once from Corollary \ref{tensor.product.several}.
\end{proof}


\begin{corollary} \label{corollary.klyachko.symmetric}
Let $\mathcal{E}$ be a toric vector bundle on the toric variety $X$. Then for each $m \in {\bf{Z}}^{+}$, the Klyachko filtrations for $Sym^{m}\mathcal{E}$ are given by
\[
\big( Sym^{m}\mathcal{E} \big)^{\rho}(i)=\sum_{i_{1}+\cdots+i_{m}=i} \operatorname{Im} \big( \mathcal{E}^{\rho}(i_{1}) \otimes \cdots \otimes \mathcal{E}^{\rho}(i_{m}) \longrightarrow Sym^{m}E \big),
\]
for each ray $\rho \in \Delta$ and each $i\in{\bf{Z}}$.
\end{corollary}
\begin{proof}
This is a particular case of the previous lemma.
\end{proof}


\begin{example} \label{example.concrete.klyachko}
Let $D = m_{1}D_{1}+\cdots+m_{d}D_{d}$ be a torus invariant Cartier divisor on the toric variety $X$, and let us denote the fiber over $t_{0}$ of the line bundle $\Ox(D)$ by $G$. The Klyachko filtration associated to a rank two toric vector bundle $\mathcal{E}$ on $X$ corresponding to a ray $\rho_{j}$ has one of the following two forms:
\begin{align*}
\begin{matrix}
{\mathcal{E}}^{\rho_{j}}(i)=
	\begin{cases}
		E& \text{if $i \leq a_{j}$},\\
		0& \text{if $i > a_{j}$},
	\end{cases} & & &
{\mathcal{E}}^{\rho_{j}}(i)=
	\begin{cases}
		E& \text{if $i \leq a_{j}$},\\
		V& \text{if $a_{j} < i \leq b_{j}$},\\
		0& \text{if $i > b_{j}$},
	\end{cases} 
\end{matrix}
\end{align*}
where $V$ is some one-dimensional subspace of the fiber $E$ of $\mathcal{E}$ over $t_{0}$, and $a_{j}$ and $b_{j}$ are some integers. For each positive integer $m$ the corresponding Klyachko filtration associated to $(Sym^{m}{\mathcal{E}}) \otimes \Ox(D)$ has respectively one of the forms: 
\begin{align*}
\begin{matrix}
{((Sym^{m}{\mathcal{E}}) \otimes \Ox(D))}^{\rho_{j}}(i)=
	\begin{cases}
		(Sym^{m}E) \otimes G & \text{if $i \leq a_{j}m+m_{j}$},\\
		0& \text{if $i > a_{j}m+m_{j}$},
	\end{cases} \\ \\
{((Sym^{m}{\mathcal{E}}) \otimes \Ox(D))}^{\rho_{j}}(i)=
	\begin{cases}
		(Sym^{m}E) \otimes G& \text{if $i \leq a_{j}m+m_{j}$},\\
		Sym^{m}_{E}\bigl(V^{\bigl\lceil \frac{i-a_{j}m-m_{j}}{b_{j}-a_{j}} \bigr\rceil}\bigr) \otimes G& \text{if $a_{j}m+m_{j} < i \leq b_{j}m+m_{j}$},\\
		0& \text{if $i > b_{j}m+m_{j}$},
	\end{cases} 
\end{matrix}
\end{align*}
where $\lceil \ \rceil$ denotes the ceiling function, and $Sym^{m}_{E}(V^{c}) =_{def} \operatorname{Im} (V^{\otimes c} \otimes E^{\otimes (m-c)} \rightarrow Sym^{m}E)$ for each integer $0 \leq c \leq m$. This convenient notation will be generalized in the next section.
\end{example}



\section{Vanishing orders on $\pe$} \label{section.vanishing.orders}

The description of the global Okounkov body of a projective variety $Z$, with respect to a flag $Y_{\bullet}$, involves identifying the image of the map $\nu_{Y_{\bullet}} \colon H^{0}(Z,\mathcal{L})\smallsetminus \{ 0 \} \longrightarrow {\bf{Z}}^{\operatorname{dim}Z}$
for each line bundle $\mathcal{L}$ on $Z$.
In this section we study these images for a suitable flag $Y_{\bullet}$, when $Z$ is the projectivization $\pe$ of a rank two toric vector bundle $\mathcal{E}$ on a smooth projective toric variety $X$. First, in \S \ref{flag} we introduce a flag of torus invariant subvarieties $Y_{\bullet}$ in $\pe$, essentially by pulling back a flag of invariant subvarieties from $X$. Next, in Definition \ref{special.sections} we present a collection of sections $\mathcal{W}_{\mathcal{L}}$ for each line bundle $\mathcal{L}$ on $\pe$. We consider these sections since we can compute their images under $\nu_{Y_{\bullet}}$ using the formulas in Lemma \ref{Lemma:OrdersOfVanishing}, and because they map onto $\nu_{Y_{\bullet}}(H^{0}(\pe,\mathcal{L})\smallsetminus \{ 0 \})$, as we prove in Proposition \ref{Proposition:OnTypesOfSectionsThatAreEnough}. In passing, we prove that after choosing an isomorphic representative of $\mathcal{L}$ so as to have an induced torus action on $H^{0}(\pe,\mathcal{L})$, the isotypical sections with respect to this action also map onto the image of $\nu_{Y_{\bullet}}$. Throughout we use the notation from \S \ref{section.notation} and additionally assume that the toric variety $X$ is smooth and projective.


\subsection{A flag of invariant subvarieties in a projectivized rank two toric vector bundle} \label{flag}
Given a toric vector bundle $\mathcal{E}$ of rank two, we construct a flag of smooth $T$-invariant subvarieties $Y_{\bullet} \colon \pe = Y_{n+1} \supseteq \ldots \supseteq Y_{0}$ in $\pe$, as follows. Let $X_{\bullet} \colon X=X_{n} \supseteq \ldots \supseteq X_{0}$ be a flag in $X$, where each $X_{i}$ is an $i$-dimensional $T$-invariant subvariety. By reordering the rays in $\Delta$ if necessary, we can assume that $X_{n-i}=\bigcap_{j=1}^{i}D_{j}$ for each $i \in \{1,\ldots,n\}$. Note that this implies that the rays $\rho_{1},\ldots,\rho_{n}$ span a maximal cone $\tau$ in $\Delta$.

Let $u_{1}$ and $u_{2}$ in $M$ be such that we can equivariantly trivialize $\mathcal{E}$ over $U_{\tau}$ as $\mathcal{E}|_{U_{\tau}} \cong \mathcal{O}_{X}(\operatorname{div}\chi^{u_{1}})|_{U_{\tau}} \oplus \mathcal{O}_{X}(\operatorname{div}\chi^{u_{2}})|_{U_{\tau}}$. The lexicographic order $\geq_{lex}$ in ${\bf{Z}}^{n}$ induces an order $\geq_{lex}$ in $M$, via the isomorphism $M \cong {\bf{Z}}^{n}$ induced by $v_{1},\ldots,v_{n}$. By reordering $u_{1}$ and $u_{2}$ if necessary, we can assume that $u_{1} \geq _{lex} u_{2}$. In other words, either $u_{1}=u_{2}$, or the first nonzero number in the list $\langle u_{1}-u_{2}, v_{1} \rangle,\ldots,\langle u_{1}-u_{2}, v_{n} \rangle$ is positive.

For each $i \in \{1,\ldots,n+1\}$, we define $Y_{i} = {\bf{P}}(\mathcal{E}|_{X_{i-1}}) = \pi^{-1}(X_{i-1}) \subseteq \pe$. To define $Y_{0}$, note that the isomorphism $\mathcal{E}|_{U_{\tau}} \cong \mathcal{O}_{X}(\operatorname{div}\chi^{u_{1}})|_{U_{\tau}} \oplus \mathcal{O}_{X}(\operatorname{div}\chi^{u_{2}})|_{U_{\tau}}$ induces an isomorphism $Y_{1} \cong {\bf{P}}(\mathcal{O}_{X_{0}} \oplus \mathcal{O}_{X_{0}})$. Hence, we get an isomorphism $\mu \colon Y_{1} \rightarrow {\bf{P}}^{1}$ between $Y_{1}$ and the projective space ${\bf{P}}^{1} \cong {\bf{P}}(\mathcal{O}_{X_{0}} \oplus \mathcal{O}_{X_{0}})$ with homogeneous coordinates $x,y$. We take $Y_{0}$ to be the point in $Y_{1}$ corresponding under $\mu$ to the point $(0:1)$, defined by the ideal $(x)$ in ${\bf{P}}^{1}$. Note that the flag $Y_{\bullet} \colon \pe = Y_{n+1} \supseteq \ldots \supseteq Y_{0}$ in $\pe$ consists of smooth $T$-invariant subvarieties.


\subsection{Computing vanishing orders} \label{subsection.vanishing.orders}
In this subsection we introduce the collection of $T$-isotypical sections $\mathcal{W}_{\mathcal{L}}$ and compute their vanishing vectors. We continue working in the setting of \S \ref{flag}.



\begin{notation} \label{notation.E.L}
We denote by $E_{1}$ and $E_{2}$ the fibers over $t_{0}$ of $\mathcal{O}_{X}(\operatorname{div}\chi^{u_{1}})$ and $\mathcal{O}_{X}(\operatorname{div}\chi^{u_{2}})$. We identify $E_{1}$ and $E_{2}$ with subspaces of the fiber $E$ of $\mathcal{E}$ over $t_{0}$ in the natural way. We denote by $L_{1},\ldots,L_{p}$ the distinct one-dimensional subspaces of $E$ that are different from $E_{1}$, but are equal to $\mathcal{E}^{\rho}(i)$ for some ray $\rho \in \Delta$ and some $i \in {\bf{Z}}$. We fix once and for all a one-dimensional subspace $L$ of $E$, different from each of the subspaces $E_{1}, L_{1}, \ldots, L_{p}$ of $E$. This is done just as an alternative to \textit{ad hoc} choices at different points in our discussion.
\end{notation}



\begin{notation} \label{notation.sym}
Let $V_{1},\ldots,V_{l}$ be subspaces of a vector space V. For any nonnegative integers $m, \linebreak[0] \alpha_{1}, \linebreak[0] \ldots, \linebreak[0] \alpha_{l}$, we define the notation $Sym_{V}^{m}(V_{1}^{\alpha_{1}} , V_{2}^{\alpha_{2}} , \ldots , V_{l}^{\alpha_{l}})$ to represent either the subspace of $Sym^{m}V$ equal to the image of the composition of the natural maps
	\[
V_{1}^{\otimes \alpha_{1}} \otimes V_{2}^{\otimes \alpha_{2}} \otimes \cdots \otimes V_{l}^{\otimes \alpha_{l}}  \otimes V^{\otimes (m - \sum_{i=1}^{l}{\alpha_{i}})} \ \longrightarrow \ V^{\otimes m} \ \longrightarrow \ Sym^{m}V,
	\]
if  $m \geq \sum_{i=1}^{l}{\alpha_{i}}$, or the subspace $0$ of $Sym^{m}V$, otherwise.
\end{notation}


On the toric variety $X$, the map defined by $(m_{n+1},\ldots,m_{d}) \mapsto \sum_{i=n+1}^{d}m_{i}D_{i}$ induces an isomorphism between ${\bf{Z}}^{d-n}$ and $\operatorname{Pic} X = N^{1}(X)$. Hence, each line bundle $\mathcal{L}$ on $\pe$ is isomorphic to a unique line bundle of the form $\Ope (m) \otimes \pi^{*} \mathcal{O}_{X}(\sum_{i=n+1}^{d}m_{i}D_{i})$.

\begin{definition} \label{special.sections} 
Let $\mathcal{L}$ be the line bundle $\mathcal{O}_{\pe}(m) \otimes \pi^{*} \mathcal{O}_{X} ( \sum_{i=n+1}^{d}{m_{i}D_{i}})$ on $\pe$, where $m, m_{n+1}, \ldots, m_{d} \in {\bf{Z}}$. Let us consider the torus action on $H^{0}\big(\pe,\mathcal{L}\big) = H^{0}\bigl(X, (Sym^{m}\mathcal{E}) \otimes  \mathcal{O}_{X}(\sum_{i=n+1}^{d}m_{i}D_{i}) \bigr)$, induced by the $T$-equivariant structure on $(Sym^{m}\mathcal{E}) \otimes  \mathcal{O}_{X}(\sum_{i=n+1}^{d}m_{i}D_{i})$. Let $G$ be the fiber over $t_{0}$ of $\Ox (\sum_{i=n+1}^{d}{m_{i}D_{i}})$. We define the following subsets of $H^{0}\big(\pe,\mathcal{L}\big)$:
\begin{align*}
\mathcal{V}_{\mathcal{L}} &= \big\{ s\in H^{0}\big(\pe,\mathcal{L}\big) \ | \ s\in H^{0}\big(\pe,\mathcal{L}\big)_{u} \smallsetminus \{ 0 \} \text{, for some $u\in M$} \big\} \\
\mathcal{W}_{\mathcal{L}} &= \big\{ s\in \mathcal{V}_{\mathcal{L}} \ | \ s(t_{0}) \text{ lies in the subspace } Sym_{E}^{m}(E_{1}^{\alpha_{0}} , L_{1}^{\alpha_{1}} , \ldots , L_{p}^{\alpha_{p}}, L^{\alpha}) \otimes G \\
	& \qquad \ \ \text{ of $Sym^{m}E$, for some $ \alpha_{0},\ldots,\alpha_{p},\alpha \in {\bf{Z}}_{\geq 0}$, satisfying $\sum_{j=0}^{p}{\alpha_{j}}+\alpha=m$} \big\}.
\end{align*}
\end{definition}


In the following lemma we give some formulas for the vanishing vector $\nu_{Y_{\bullet}}(s)$ of a section $s \in \mathcal{W}_{\mathcal{L}}$. 

\begin{lemma}[]\label{Lemma:OrdersOfVanishing}
Let $\mathcal{L}$ be the line bundle $\Ope(m) \otimes \pi^{*}\Ox(\sum_{i=n+1}^{d}{m_{i}D_{i}})$ on $\pe$,
for some $m, m_{n+1},\ldots,m_{d} \in {\bf{Z}}$. Let $G$ be the fiber over $t_{0}$ of $\Ox (\sum_{i=n+1}^{d}{m_{i}D_{i}})$. Let $s$ be a nonzero section in $H^{0} \big( \pe, \mathcal{L} \big)_{u}=H^{0} \big( X, ( Sym^{m}\mathcal{E} ) \otimes \Ox (\sum_{i=n+1}^{d}{m_{i}D_{i}}) \big)_{u}$, for some $u \in M$, and let $\nu_{Y_{\bullet}}(s) = (\nu_{1}, \ldots , \nu_{n+1}) \in {\bf{Z}}^{n+1} $. Then: 
			\item[\textbf{(a)}] For each $j \in \{ 1,\ldots,n \}$ we have
$
\nu_{j}=\langle \nu_{n+1}u_{1}+(m-\nu_{n+1})u_{2}-u , v_{j} \rangle.
$
			\item[\textbf{(b)}] If $s(t_{0})$ lies in the subspace $Sym_{E}^{m}(E_{1}^{\alpha_{0}} , L_{1}^{\alpha_{1}} , \ldots , L_{p}^{\alpha_{p}}, L^{\alpha}) \otimes G$ of $(Sym^{m} E) \otimes G$, for some $\alpha_{0},\ldots,\alpha_{p},\alpha \in {\bf{Z}}_{\geq 0}$ such that $\sum_{i=0}^{p}{\alpha_{i}}+\alpha=m$, then $ \nu_{n+1} = \alpha_{0}$.		
\end{lemma}
\begin{proof}
			\item[\textbf{(a)}] The vector $\nu_{Y_{\bullet}}(s)$ can be computed in any neighborhood of $Y_{0}$ in $\pe$. Hence we can assume that $X=U_{\tau}$, that $\mathcal{E}=\Ox(\operatorname{div}\chi^{u_{1}}) \oplus \Ox(\operatorname{div}\chi^{u_{2}}) $, and that $s$ is a section in $H^{0} \big( X, Sym^{m}\mathcal{E}\big)_{u}$. Note that $Sym^{m}\mathcal{E}= \bigoplus^{m}_{i=0} \Ox(\operatorname{div}\chi^{(m-i)u_{1}+iu_{2}})$, and so $s$ corresponds to the section 
\begin{equation*}
\big{(} c_{0}\chi^{-u} ,\ldots, c_{i}\chi^{-u} ,\ldots, c_{m}\chi^{-u} \big{)} \in \bigoplus^{m}_{i=0} H^{0} \big(X,\Ox(\operatorname{div}\chi^{(m-i)u_{1}+iu_{2}}) \big)_{u} =H^{0} \big( X, Sym^{m}\mathcal{E}\big)_{u},
\end{equation*}
for some $c_{0},\ldots,c_{m}\in k$. Let us denote $\Ox \oplus \Ox$ by $\mathcal{E}'$. By combining the natural isomorphisms in each component, we get an isomorphism $\mathcal{E}=\Ox(\operatorname{div}\chi^{u_{1}}) \oplus \Ox(\operatorname{div}\chi^{u_{2}}) \cong \Ox \oplus \Ox = \mathcal{E}'$. This isomorphism induces a commutative diagram,
\[
\xymatrix{
\pe = {\bf{P}} \big(\Ox(\operatorname{div}\chi^{u_{1}}) \oplus \Ox(\operatorname{div}\chi^{u_{2}}) \big) \ar@/_/[dr] \ar[rr]^-{\varphi} & &{\bf{P}} \big( \Ox \oplus \Ox \big) = {\bf{P}}(\mathcal{E}') \ar@/^/[ld] \\
& X  &}
\]
where the map $\varphi$ is an isomorphism. Let $Y_{\bullet}' \colon Y'_{n+1} \supseteq Y'_{n} \supseteq \cdots \supseteq Y'_{0}$ be the $T$-invariant flag in ${\bf{P}} \big( \mathcal{E}' \big)$, as defined in \S \ref{flag}. Note that the $T$-invariant flags in $\pe$ and ${\bf{P}} ( \mathcal{E}' )$ correspond to each other under $\varphi$. On $Y'_{i} = {\bf{P}}(\mathcal{E}'|_{X_{i-1}})$, let us denote $\mathcal{O}_{{\bf{P}}(\mathcal{E}'|_{X_{i-1}})}(m)$ by $\mathcal{O}_{Y'_{i}}(m)$ for each $i \in \{ 1,\ldots,n+1 \}$. Under the isomorphism $\varphi$, $s$ corresponds to the section 
\[
\begin{split}
s' = s'_{n+1} &= \big{(} c_{0}\chi^{mu_{1}-u} ,\ldots, c_{i}\chi^{(m-i)u_{1}+iu_{2}-u} ,\ldots, c_{m}\chi^{mu_{2}-u} \big{)} \in \bigoplus^{m}_{i=0} H^{0} \big(X, \mathcal{O}_{X} \big) \\
&= H^{0} \big( X, Sym^{m}(\mathcal{E}')\big) = H^{0} \big ({\bf{P}}(\mathcal{E}'),\mathcal{O}_{{\bf{P}}(\mathcal{E}')}(m) \big) = H^{0} \big( Y'_{n+1},\mathcal{O}_{Y'_{n+1}} (m) \big). 
\end{split}
\]
Note that $(\nu_{1},\ldots,\nu_{n+1}) = \nu_{Y_{\bullet}}(s)= \nu_{Y'_{\bullet}}(s')$. Let us set $h=\operatorname{max} \, \{i \ | \ 0 \leq i \leq m \text{ and } c_{i} \neq 0 \}$, and let $v_{1}^{*},\ldots,v_{n}^{*}$ be the basis of $M$ dual to the basis $v_{1},\ldots,v_{n}$ of $N$. 

It is straightforward to see that when we follow the procedure to compute $\nu_{Y'_{\bullet}}(s')$ outlined in \S \ref{okounkov.procedure}, for each $0 \leq l \leq n$, the section obtained in the step when we restrict to $Y'_{n+1-l}$ corresponds to the section
\begin{equation} \label{s'}
\begin{split}
s'_{n+1-l} &= \big{(} c_{0}\chi^{mu_{1}-u-\sum_{j=1}^{l}{\nu_{j}v^{*}_{j}}}|_{X_{n-l}} ,\ldots,  c_{i}\chi^{(m-i)u_{1}+iu_{2}-u-\sum_{j=1}^{l}{\nu_{j}v^{*}_{j}}}|_{X_{n-l}} ,\ldots, \\ & \qquad c_{m}\chi^{mu_{2}-u-\sum_{j=1}^{l}{\nu_{j}v^{*}_{j}}}|_{X_{n-l}} \big{)} \in \bigoplus^{m}_{i=0} H^{0} \big(X_{n-l},  \mathcal{O}_{X_{n-l}} \big) = H^{0} \big( X_{n-l}, Sym^{m}(\mathcal{E}'|_{X_{n-l}} )\big) \\ & \qquad \qquad = H^{0} \big ({\bf{P}}(\mathcal{E}'|_{X_{n-l}}),\mathcal{O}_{{\bf{P}}(\mathcal{E}'|_{X_{n-l}})}(m) \big) = H^{0} \big( Y'_{n+1-l},\mathcal{O}_{Y'_{n+1-l}} (m) \big), 
\end{split}
\end{equation}
under the natural identification.

Assume now that for some $0 \leq l < n$, we have proven that $\nu_{j}=\langle (m-h)u_{1}+hu_{2}-u , v_{j} \rangle$, for each $j \in \{ 1,\ldots,l \}$.
Note that for each $a \in {\bf{Z}}_{\geq 0}$ we have the following commutative diagram
\[
\begin{CD}
H^{0} \big( Y'_{n+1-l}, \mathcal{O}_{Y'_{n+1-l}}(m) \otimes \mathcal{O}_{Y'_{n+1-l}} (-aY'_{n-l}) \big) @>\varphi_{a}>> H^{0} \big( Y_{n+1-l}, \mathcal{O}_{Y'_{n+1-l}}(m) \big) \\
@| @| \\
H^{0} \big( X_{n-l}, Sym^{m}(\mathcal{E}'|_{X_{n-l}}) \otimes \mathcal{O}_{X_{n-l}} (-aX_{n-l-1}) \big) @>>> H^{0} \big( X_{n-l}, Sym^{m}(\mathcal{E}'|_{X_{n-l}}) \big),
\end{CD}
\]
and denote the map in its top row by $\varphi_{a}$.
Next, we note that
\begin{equation}
\begin{split} \label{gamma}
\nu_{l+1} &= \operatorname{max} \, \{a \in {\bf{Z}}_{\geq 0} \ | \ s'_{n+1-l} \in \operatorname{Im}(\varphi_{a}) \}  \\
&= \operatorname{max} \, \{a \in {\bf{Z}}_{\geq 0} \ | \ c_{i} \chi^{(m-i)u_{1}+iu_{2}-u- \sum_{j=1}^{l}{\nu_{j}v^{*}_{j}}}|_{X_{n-l}} \in \operatorname{Im}(H^{0}(X_{n-l}, \mathcal{O}_{X_{n-l}} (-aX_{n-l-1})) \hookrightarrow  \\ 
& \phantom{= \operatorname{max} \, \{a \in {\bf{Z}}_{\geq 0} \ | \ } H^{0}(X_{n-l}, \mathcal{O}_{X_{n-l}} )) \text{ for each $i=0,\ldots,m$}\} \\
&= \operatorname{max} \, \{a \in {\bf{Z}}_{\geq 0} \ | \ a \leq \langle (m-i)u_{1}+iu_{2}-u , v_{l+1} \rangle \text{ for each } i \text { such that } \\
& \phantom{= \operatorname{max} \, \{a \in {\bf{Z}}_{\geq 0} \ | \ \ } c_{i} \chi^{(m-i)u_{1}+iu_{2}-u- \sum_{j=1}^{l}{\nu_{j}v^{*}_{j}}}|_{X_{n-l}}  \neq 0 \}.
\end{split}
\end{equation}
We also note that $c_{h} \chi^{(m-h)u_{1}+hu_{2}-u- \sum_{j=1}^{l}{\nu_{j}v^{*}_{j}}}|_{X_{n-l}} \neq 0$. Now, if $\langle u_{1} - u_{2} , v_{l+1} \rangle < 0$, then there exists $q \in \{1, \ldots, l \}$, such that $\langle u_{1} - u_{2} , v_{q} \rangle > 0$. In this case, for each $i \in \{0,\ldots, h-1 \}$, it follows that $c_{i} \chi^{(m-i)u_{1}+ iu_{2}- u - \sum_{j=1}^{l}{\nu_{j}v^{*}_{j}}}|_{X_{n-q}} = 0$. Hence, either $\langle u_{1} - u_{2} , v_{l+1} \rangle \geq 0$, or $c_{i} \chi^{(m-i)u_{1}+ iu_{2}- u - \sum_{j=1}^{l}{\nu_{j}v^{*}_{j}}}|_{X_{n-l}} = 0$ for each $i \in \{ 0,\ldots, h-1 \}$. In either case, it follows from (\ref{gamma}) that 
$
\nu_{l+1} = \langle (m-h)u_{1}+hu_{2}-u , v_{l+1} \rangle.
$
Therefore we can iterate this procedure, and in this way we obtain that for each $j \in \{ 1,\ldots,n \}$, 
\[
\nu_{j}=\langle (m-h)u_{1}+hu_{2}-u , v_{j} \rangle.
\]
Now, $\nu_{n+1}$ is equal to the vanishing order along $Y'_{0}$ of the section $s'_{1} \in H^{0}(Y'_{1},\mathcal{O}_{Y'_{1}}(m))$ described in (\ref{s'}) for $l=n$.
We have a natural isomorphism $\mu \colon Y'_{1} \rightarrow {\bf{P}}^{1}$ between $Y'_{1} = {\bf{P}}(\mathcal{O}_{X_{0}} \oplus \mathcal{O}_{X_{0}})$ and the projective space ${\bf{P}}^{1} \cong {\bf{P}}(\mathcal{O}_{X_{0}} \oplus \mathcal{O}_{X_{0}})$ with homogeneous coordinates $x,y$. Recall that under $\mu$, $Y'_{0}$ corresponds to the point $(0:1)$, defined by the ideal $(x)$ in ${\bf{P}}^{1}$. Under $\mu$, $\mathcal{O}_{Y'_{1}}(m)$ corresponds to $\mathcal{O}_{{\bf{P}}^1}(m)$. Depending on whether or not $u_{1}=u_{2}$, there are two possibilities for the section in $H^{0}({\bf{P}}^{1},\mathcal{O}_{{\bf{P}}^{1}}(m))$ that corresponds to $s'_{1}$. Namely, $s'_{1}$ corresponds to $c_{h}x^{m-h}y^{h}$ if $u_{1} \neq u_{2}$, and it corresponds to $\sum^{m}_{i=0}{c_{i}x^{m-i}y^{i}}$ if $u_{1} = u_{2}$. In either case, we obtain that $\nu_{n+1}=m-h$, and then part (a) is proven.
			\item[\textbf{(b)}] As in part (a), we can reduce to  the case when $X=U_{\tau}$, $\mathcal{E}=\Ox(\operatorname{div}\chi^{u_{1}}) \oplus \Ox(\operatorname{div}\chi^{u_{2}}) $, $s \in H^{0} \big( X, Sym^{m}\mathcal{E}\big)_{u}$, and $s(t_{0})$ lies in the subspace $Sym_{E}^{m}(E_{1}^{\alpha_{0}} , L_{1}^{\alpha_{1}} , \ldots , L_{p}^{\alpha_{p}}, L^{\alpha})$ of $Sym^{m}E$.
Let $x,y \in E$ be such that $E_{1} = kx$ and $E_{2}= ky$. Then $x$ and $y$ form a basis for $E$, and $x^{m-i}y^{i}$ for $i=0,\ldots,m$ form a basis for $Sym^{m}E$. Let $\beta_{1},\ldots,\beta_{p},\beta \in k$ be such that $L = k(\beta x+y)$ and $L_{i} = k(\beta_{i}x+y)$ for each $i \in \{ 1,\ldots,p \}$.
For $c_{0},\ldots,c_{m}\in k$, defined as in part (a), we proved that $\operatorname{max} \, \{i \ | \ 0 \leq i \leq m \text{ and } c_{i} \neq 0 \} = m - \nu_{n+1}$. On the one hand, we see that the image of $s$ at $t_{0}$ is $s(t_{0}) = \sum_{i=0}^{m}{c_{i}x^{m-i}y^{i}} \in Sym^{m}E$. On the other hand, 
\[
s(t_{0}) \in Sym_{E}^{m}(E_{1}^{\alpha_{0}} , L_{1}^{\alpha_{1}} , \ldots , L_{p}^{\alpha_{p}}, L^{\alpha}) = k(\sum_{ i =0 }^{m-\alpha_{0}-1}{\beta'_{i}}x^{m-i}y^{i}+x^{\alpha_{0}}y^{m-\alpha_{0}}),
\]
for some $\beta'_{0},\ldots,\beta'_{m-\alpha_{0}-1} \in k$. From this it follows that $\nu_{n+1} = \alpha_{0}$ as desired.
\end{proof}


\subsection{The image of $\nu_{Y_{\bullet}}$} \label{subsection.enough.sections}

In this subsection we prove that $\nu_{Y_{\bullet}}$ maps the collection of $T$-isotypical sections $\mathcal{W}_{\mathcal{L}}$ onto $\nu_{Y_{\bullet}}(H^{0}(\pe,\mathcal{L})\smallsetminus \{ 0 \})$. We continue working in the setting of \S \ref{flag}-\ref{subsection.vanishing.orders}.

We start by proving that if a line bundle $\mathcal{L}$ and a flag $Y_{\bullet}$ on an affine variety $Z$ are suitably compatible with the action of a torus $T$ on $Z$, then the nonzero $T$-isotypical sections of $\mathcal{L}$ map onto the image of $\nu_{Y_{\bullet}}$.

\begin{lemma}[] \label{preLemma:EquivariantSectionsAreEnough}
Let $Z$ be an affine variety with an algebraic action of a torus $T$, and a flag $Y_{\bullet} \colon Z=Y_{l} \supseteq Y_{l-1} \supseteq \ldots \supseteq Y_{0}$, where each $Y_{i}$ is a normal $i$-dimensional $T$-invariant subvariety.
Assume that for each $i \in \{ 1,\ldots,l \}$, there is a $T$-isotypical rational function $h_{i}$ on $Y_{i}$,
such that $Y_{i-1} = \operatorname{div}h_{i}$. Let $g$ be a $T$-isotypical rational function on $Z$ and let $s_{1},\ldots,s_{q} \in H^{0}(Z,\mathcal{O}_{Z}(\operatorname{div}g))$ be nonzero $T$-isotypical sections corresponding to distinct characters of $T$. Then
$\nu_{Y_{\bullet}}(s_{1}+\cdots+s_{q}) \in \{\nu_{Y_{\bullet}}(s_{1}),\ldots, \nu_{Y_{\bullet}}(s_{q})\}$.
\end{lemma}
\begin{proof}
We proceed by induction on the dimension of $Z$. Let $Z_{\bullet} \colon Y_{l-1}=Z_{l-1} \supseteq Z_{l-2} \supseteq \ldots \supseteq Z_{0}$ be the flag of normal $T$-invariant subvarieties in $Y_{l-1}$ defined by $Z_{i} = Y_{i}$, for each $i$. Using the natural isomorphism $H^{0}(Z,{O}_{Z}(\operatorname{div}g)) \cong H^{0}(Z,\mathcal{O}_{Z})$, we can reduce to the case when $g=1$ and the sections are identified with regular functions. Let $s=s_{1}+\cdots+s_{q}$. For each $a \in {\bf{Z}}_{\geq 0}$ the natural inclusion map
$
\varphi_{a} \colon H^{0}(Y_{l},\mathcal{O}_{Y_{l}}(-aY_{l-1})) \longhookrightarrow H^{0}(Y_{l},\mathcal{O}_{Z})
$
is compatible with the decomposition of these spaces into $T$-isotypical summands. It follows that 
\[ 
\nu_{1}(s) = \operatorname{ord}_{Y_{l-1}}(s) = \operatorname{min}  \{\operatorname{ord}_{Y_{l-1}}(s_{1}),\ldots, \operatorname{ord}_{Y_{l-1}}(s_{q})\} = \operatorname{min}  \{\nu_{1}(s_{1}),\ldots, \nu_{1}(s_{q})\}.
\]
If we reorder the sections so that $\nu_{1}(s) = \nu_{1}(s_{i})$ for $1 \leq i \leq e$, and $\nu_{1}(s) < \nu_{1}(s_{i})$ for $e+1 \leq i \leq q$, for some $e \in \{1,\ldots, q \}$, then $(h_{l}^{-\nu_{1}(s)}s)|_{Y_{l-1}} = \sum_{i=1}^{e} (h_{l}^{-\nu_{1}(s)}s_{i})|_{Y_{l-1}}$ and using the induction hypothesis we get
\begin{align*}
\nu_{Y_{\bullet}}(s) = \bigl(\nu_{1}(s),\nu_{Z_{\bullet}} \bigl((h_{l}^{-\nu_{1}(s)}s)|_{Y_{l-1}} \bigr)\bigr) &\in \bigl\{
\bigl(\nu_{1}(s),\nu_{Z_{\bullet}} \bigl((h_{l}^{-\nu_{1}(s)}s_{i})|_{Y_{l-1}} \bigr)\bigr) \ | \ i =1,\ldots,e \bigr\} \\
& = \{ \nu_{Y_{\bullet}} (s_{1}) , \ldots, \nu_{Y_{\bullet}} (s_{e}) \} \subseteq \{ \nu_{Y_{\bullet}} (s_{1}) , \ldots, \nu_{Y_{\bullet}} (s_{q}) \}.
\end{align*}
\end{proof}


Recall that any line bundle on $\pe$ is isomorphic to a unique line bundle of the form $\mathcal{O}_{\pe}(m) \otimes \pi^{*} \mathcal{O}_{X}(\sum_{i=n+1}^{d}{m_{i}D_{i}})$. In the following proposition, we prove that for a line bundle on $\pe$ of that form, the $T$-isotypical sections map onto the image of $\nu_{Y_{\bullet}}$. 

\begin{proposition} \label{Lemma:EquivariantSectionsAreEnough}
Let $\mathcal{L}$ be the line bundle $\mathcal{O}_{\pe}(m) \otimes \pi^{*} \mathcal{O}_{X}(\sum_{i=n+1}^{d}{m_{i}D_{i}})$ on $\pe$ for some $m, m_{1}, \ldots, m_{d} \in {\bf{Z}}$, and let $s$ be a nonzero global section of $\mathcal{L}$. Let $s_{1},\ldots,s_{q} \in H^{0}(\pe,\mathcal{L})$ be the unique nonzero $T$-isotypical sections corresponding to distinct characters of $T$ such that $s=s_{1}+\cdots+s_{q}$. Then $\nu_{Y_{\bullet}}(s) \in \{\nu_{Y_{\bullet}}(s_{1}),\ldots, \nu_{Y_{\bullet}}(s_{q})\}$.
\end{proposition}
\begin{proof}
It is enough to consider the case $X = U_{\tau}$. There is a natural choice of coordinates $X = \operatorname{Spec}k[x_{1},\ldots,x_{n}]= {\bf{A}}^{n}$ and ${\bf{P}}(\mathcal{E}) = \operatorname{Spec}k[x_{1},\ldots,x_{n}] \times \operatorname{Proj}k[x,y]= {\bf{A}}^{n} \times {\bf{P}}^{1}$, which is induced by the ordering of the rays of $\tau$ and the trivialization of $\mathcal{E}$ over $U_{\tau}$. In these coordinates we have that $X_{i}  = \{ (x_{1},\ldots,x_{n}) \in {\bf{A}}^{n} \ | \ x_{j} = 0 \textnormal{ for } 1\leq j \leq n-i   \}$, and $Y_{i+1} = X_{i} \times {\bf{P}}^{1}$ for $0 \leq i \leq n$ and $Y_{0} = X_{0} \times \{ (0:1) \}$. We also have that $T = \operatorname{Spec}k[x_{1},\ldots,x_{n}]_{x_{1} \cdots x_{n}} = (k^{*})^{n}$ acts on ${\bf{A}}^{n}$ by componentwise multiplication, and that an element $t=(t_{1},\ldots,t_{n}) \in T$ acts on $P = \bigl( (x_{1},\ldots,x_{n}) , (x:y) \bigr) \in {\bf{A}}^{n} \times {\bf{P}}^{1}$ by 
$
tP= \bigl( (t_{1}x_{1},\ldots,t_{n}x_{n}) , (t_{1}^{ \langle u_{1}, v_{1} \rangle } \cdots t_{n}^{ \langle u_{1}, v_{n} \rangle } x : t_{1}^{ \langle u_{2}, v_{1} \rangle } \cdots t_{n}^{ \langle u_{2}, v_{n} \rangle }y) \bigr). 
$
If $U \subseteq {\bf{P}}^{1}$ is the complement of $(1:0)$, it is enough to prove that in the $T$-invariant affine open set $Z = {\bf{A}}^{n} \times U$ the restriction of the flag $Y_{\bullet}$ and the line bundle $\mathcal{L}$ satisfy the hypotheses of Lemma \ref{preLemma:EquivariantSectionsAreEnough}. 
We show that $\Ope(1)|_{Z} = \mathcal{O}_{Z} (\operatorname{div}g)$ for some $T$-isotypical rational function $g$ on $Z$, since from this all the assertions follow at once. 
The surjection $\mathcal{E} \rightarrow \mathcal{O}_{X}(\operatorname{div}\chi^{u_{2}})$ corresponds to a geometric section $\theta$ of the projection $\pi \colon \pe \rightarrow X$, i.e. a morphism $\theta \colon X \rightarrow \pe$ such that $\pi \circ \theta =id_{X}$. If we set $X' = \theta(X)$, then $\pi_{*}(\mathcal{O}_{\pe}(1) \otimes \mathcal{O}_{X'}) = \mathcal{O}_{X}(\operatorname{div} \chi^{u_{2}})$ and we have the exact sequence
\begin{equation} \label{sequence}
0 \longrightarrow \mathcal{O}_{\pe}(1) \otimes \mathcal{O}_{\pe}(-X') \longrightarrow \mathcal{O}_{\pe}(1) \longrightarrow \mathcal{O}_{\pe}(1) \otimes \mathcal{O}_{X'} \longrightarrow 0. 
\end{equation}
By Grauert's theorem (see III.12.9 in \cite{Hartshorne}) $R^{1}\pi_{*} (\mathcal{O}_{\pe}(1) \otimes \mathcal{O}_{\pe}(-X'))=0$. Then applying $\pi_{*}$ to (\ref{sequence}) gives $\pi_{*}(\mathcal{O}_{\pe}(1) \otimes \mathcal{O}_{\pe}(-X')) = \mathcal{O}_{X}(\operatorname{div}\chi^{u_{1}})$. 
Since $\mathcal{O}_{\pe}(1) \otimes \mathcal{O}_{\pe}(-X')$ has degree zero along the fibers of $\pi$, 
\[
\mathcal{O}_{\pe}(1) \otimes \mathcal{O}_{\pe}(-X') = \pi^{*}\pi_{*}(\mathcal{O}_{\pe}(1) \otimes \mathcal{O}_{\pe}(-X'))
 = \pi^{*} \mathcal{O}_{X}(\operatorname{div}\chi^{u_{1}}). 
\]
And since in local coordinates $X' = {\bf{A}}^{n} \times (0:1)$, we can take $g =  (x/y) \pi^{*}(\chi^{u_{1}})$.
\end{proof}


The following lemma will be used in the proofs of Proposition \ref{Proposition:OnTypesOfSectionsThatAreEnough} and Theorem \ref{theorem}.

\begin{lemma} \label{Lemma:SubspacesOfSymE}
Let $V_{1},\ldots,V_{l}$ be distinct one-dimensional subspaces of a two-dimensional vector space $V$. Let $m, \alpha_{1}, \ldots, \alpha_{l}$ be nonnegative integers. Then 
\[
Sym_{V}^{m}(V_{1}^{\alpha_{1}}) \cap Sym_{V}^{m}(V_{2}^{\alpha_{2}}) \cap \cdots \cap Sym_{V}^{m}(V_{l}^{\alpha_{l}}) = Sym_{V}^{m}(V_{1}^{\alpha_{1}} , V_{2}^{\alpha_{2}} , \ldots , V_{l}^{\alpha_{l}}).
\]  
Furthermore, this subspace of $Sym^{m}V$ is nonzero precisely when $m \geq \sum^{l}_{i=1}{\alpha_{i}}$, and in that case its dimension is $m+1-\sum^{l}_{i=1}{\alpha_{i}}$.
\begin{proof}
We fix an isomorphism of $k$-algebras between $\oplus_{h\geq 0} Sym^{h}V$ and the polynomial ring in two variables $k[x,y]$. The subspaces $V_{1},\ldots,V_{l}$ of $V$ correspond to the linear spans of some distinct linear forms $f_{1},\ldots,f_{l}$. The subspaces $ \cap_{i=1}^{l}Sym_{V}^{m}(V_{i}^{\alpha_{i}})$ and $Sym_{V}^{m}(V_{1}^{\alpha_{1}} , \ldots , V_{l}^{\alpha_{l}})$ of $Sym^{m}V$ both correspond to the homogeneous polynomials of degree $m$ divisible by $f_{1}^{\alpha_{1}}\cdots f_{l}^{\alpha_{l}}$. From this observation the conclusion follows.
\end{proof} 
\end{lemma}


In the next proposition we prove that for every line bundle $\mathcal{L}$ on $\pe$, in order to find the image of $\nu_{Y_{\bullet}}$, we can restrict our attention to the sections in $\mathcal{W}_{\mathcal{L}}$.

\begin{proposition}[]\label{Proposition:OnTypesOfSectionsThatAreEnough}
Let $\mathcal{L}$ be the line bundle $\mathcal{L}=\Ope(m) \otimes \pi^{*} \Ox (\sum_{i=n+1}^{d}{m_{i}D_{i}})$ on $\pe$, for some $m, m_{n+1},\ldots,m_{d} \in {\bf{Z}}$. Then we have the following equality of subsets of ${\bf{Z}}^{n+1}$:
\[ 
\nu_{Y_{\bullet}}(H^{0}\big(\pe,\mathcal{L}\big) \smallsetminus \{ 0 \}) = \nu_{Y_{\bullet}}(\mathcal{V}_{\mathcal{L}}) = \nu_{Y_{\bullet}}(\mathcal{W}_{\mathcal{L}}).
\]
\begin{proof}
From their definitions, we have that $\nu_{Y_{\bullet}}(H^{0}\big(\pe,\mathcal{L}\big) \smallsetminus \{ 0 \}) \supseteq \nu_{Y_{\bullet}}(\mathcal{V}_{\mathcal{L}}) \supseteq \nu_{Y_{\bullet}}(\mathcal{W}_{\mathcal{L}})$. The sets $\nu_{Y_{\bullet}}(H^{0}\big(\pe,\mathcal{L}\big) \smallsetminus \{ 0 \})$ and $\nu_{Y_{\bullet}}(\mathcal{V}_{\mathcal{L}})$ are equal by Proposition \ref{Lemma:EquivariantSectionsAreEnough}. Let us consider $\nu_{Y_{\bullet}}(s) \in \nu_{Y_{\bullet}}(H^{0}\big(\pe,\mathcal{L}\big) \smallsetminus \{ 0 \}) = \nu_{Y_{\bullet}}(\mathcal{V}_{\mathcal{L}})$. We can assume that $s \in H^{0}\big(\pe,\mathcal{L} \big)_{u}$, for some $u\in M$. By Remark \ref{Lemma:LazMus}, the set $\nu_{Y_{\bullet}}(H^{0}(\pe,\mathcal{L})_{u} \smallsetminus \{ 0 \})$ is finite with cardinality $\operatorname{dim}_{k}H^{0}\big(\pe,\mathcal{L} \big)_{u}$.  Let us denote the fiber over $t_{0}$ of $\Ox (\sum_{i=n+1}^{d} {m_{i}D_{i}})$ by $G$. From Example  \ref{example.concrete.klyachko} and Lemma \ref{Lemma:SubspacesOfSymE} we see that there exist $\alpha_{0},$ $\ldots,$ $\alpha_{p}$ $\in$ ${\bf{Z}}_{\geq 0}$ such that 
\begin{align*}
 \operatorname{Im} \Big( H^{0}\big(\pe,\mathcal{L} \big)_{u} &= H^{0} \big( X,\pi_{*}\mathcal{L} \big)_{u} \longhookrightarrow (Sym^{m}E) \otimes G \Big) \\ 
&=  (\pi_{*}\mathcal{L})^{\rho_{1}}(\langle u, v_{1}  \rangle ) \cap (\pi_{*}\mathcal{L})^{\rho_{2}}(\langle u, v_{2}  \rangle ) \cap \cdots \cap (\pi_{*}\mathcal{L})^{\rho_{d}}(\langle u, v_{d}  \rangle )  \\
& = (Sym_{E}^{m}(E_{1}^{\alpha_{0}}) \otimes G) \cap (Sym_{E}^{m}(L_{1}^{\alpha_{1}}) \otimes G) \cap \cdots \cap (Sym_{E}^{m}(L_{p}^{\alpha_{p}}) \otimes G) \\
& = Sym_{E}^{m}(E_{1}^{\alpha_{0}},L_{1}^{\alpha_{1}},L_{2}^{\alpha_{2}},\ldots,L_{p}^{\alpha_{p}}) \otimes G.
\end{align*}
Using again Lemma \ref{Lemma:SubspacesOfSymE} we see that $\alpha =_{def} m-\sum_{i=0}^{p}{\alpha_{i}}$ is a non negative integer and $\alpha +1 = \operatorname{dim}_{k}H^{0}\big(\pe,\mathcal{L} \big)_{u} $. For each $j \in \{0,\ldots,\alpha \}$, let $s_{j} \in    H^{0}\big(X,\pi_{*}\mathcal{L} \big)_{u} \smallsetminus \{ 0 \} = H^{0}\big(\pe,\mathcal{L} \big)_{u} \smallsetminus \{ 0 \}$ be such that 
\[
s_{j}(t_{0}) \in Sym_{E}^{m}(E_{1}^{\alpha_{0}+j},L_{1}^{\alpha_{1}},L_{2}^{\alpha_{2}},\ldots,L_{p}^{\alpha_{p}},L^{\alpha-j}) \otimes G.
\]
From Lemma \ref{Lemma:OrdersOfVanishing} it follows that $\nu_{Y_{\bullet}}(s_{0}),\ldots,\nu_{Y_{\bullet}}(s_{\alpha})\in \nu_{Y_{\bullet}}(H^{0}\big(\pe,\mathcal{L}\big)_{u} \smallsetminus \{ 0 \})$ are pairwise distinct, so
\[
\nu_{Y_{\bullet}}(s) \in \nu_{Y_{\bullet}}(H^{0}\big(\pe,\mathcal{L}\big)_{u} \smallsetminus \{ 0 \}) = \{ \nu_{Y_{\bullet}}(s_{0}),\ldots,\nu_{Y_{\bullet}}(s_{\alpha}) \} \subseteq \nu_{Y_{\bullet}}(\mathcal{W}_{\mathcal{L}}),
\] 
and this completes the proof of the proposition.
\end{proof}
\end{proposition}



\section{The global Okounkov body of $\pe$} \label{section.okounkov.body}

In this section we describe the global Okounkov body of the projectivization of a rank two toric vector bundle over a smooth projective toric variety, with respect to the flag of invariant subvarieties constructed in \S \ref{flag}. We introduce the relevant terminology in \S \ref{inequalities}, and we prove our result describing this global Okounkov body in terms of linear inequalities in \S \ref{subsection.theorem}. Throughout this section we use the notation and constructions from \S \ref{section.notation} and \S \ref{flag}-\ref{subsection.vanishing.orders}.


\subsection{Supporting hyperplanes of the global Okounkov body of $\pe$} \label{inequalities}
Let $\mathcal{E}$ be a toric vector bundle of rank two over the smooth projective toric variety $X$. Let $u_{1}, u_{2} \in M$, with $u_{1} \geq _{lex} u_{2}$, be as defined in \S \ref{flag}, and let the subspaces $E_{1},L_{1},\ldots,L_{p}$ of the fiber $E$ of $\mathcal{E}$ over $t_{0}$ be as defined in \S \ref{subsection.vanishing.orders}. Let us classify the filtrations $\{\mathcal{E}^{\rho_{j}} \ | \ j=1,\ldots,d   \}$ associated to $\mathcal{E}$ by defining
\begin{align*}
A &=_{def}\{j \in \{1,\ldots,d\} \ | \ \operatorname{dim}_{k}{\mathcal{E}}^{\rho_{j}}(i) \neq 1 \text{ for all } i\in {\bf{Z}} \} \\
B &=_{def}\{j \in \{1,\ldots,d\} \ | \ \mathcal{E}^{\rho_{j}}(i) = E_{1} \text{ for some } i\in {\bf{Z}} \} \\
C_{h} &=_{def}\{ j \in \{ 1,\ldots,d \} \ | \ \mathcal{E}^{\rho_{j}}(i)= L_{h} \textnormal{ for some $i \in {\bf{Z}}$}  \}
\end{align*}
for each $h \in \{ 1,\ldots , p \}$. And let us define a nonempty set $J \subseteq \{1,\ldots,d\}$ to be \emph{admissible} if it has one of the following three forms:
\begin{itemize}
	\item $J=\{j\}$ for some $j \in A$.
\sbl	
	\item $J=\{j\}$ for some $j \in B$.
\sbl
	\item $J=\{j_{1},\ldots,j_{l}\}$ for some $j_{1},\ldots,j_{l} \in \{1,\ldots,d\}$ such that there exist distinct indices $i_{1},\ldots,i_{l} \in \{1,\ldots,p\}$ with $j_{h} \in C_{i_{h}}$, for each $h \in \{ 1,\ldots,l \}$.
\end{itemize}
Note that each admissible subset of $\{1,\ldots,d\}$ is contained in exactly one of the sets $A$, $B$ and $C =_{def}\cup_{i=1}^{p}C_{i}$.
For each ray $\rho_{j} \in \Delta$ we define integers $a_{j}$ and $b_{j}$ as follows. Let $a_{j}=\operatorname{max}\{i\in {\bf{Z}} \ | \ \mathcal{E}^{\rho_{j}}(i)=E \}$, and let
\begin{equation*} 
b_{j}=
	\begin{cases}
		a_{j}+1 & \text{if $j \in A$}, \\
		\operatorname{max}\{i\in {\bf{Z}} \ | \ \operatorname{dim}_{k}{\mathcal{E}}^{\rho_{j}}(i) = 1 \} & \text{if $j \in B \cup C$}.
	\end{cases}
\end{equation*}

\begin{example} \label{tangent.bundle.definition.2}
In the case of $T_{{\bf{P}}^{2}}$, the tangent bundle of the projective plane (see Example \ref{tangent.bundle.definition}), if we take the $T$-invariant flag ${\bf{P}}^{2} \supseteq D_{1} \supseteq D_{1} \cap D_{2}$ in ${\bf{P}}^{2}$, we get $\tau = \sigma_{3}$, $u_{1} = (1,0)$ and $u_{2} = (0,1)$. We get $a_{1}=a_{2}=a_{3}=0$ and $b_{1}=b_{2}=b_{3}=1$. We also get $E_{1} = V_{1}$, $L_{1} = V_{2}$ and $L_{2} = V_{3}$, and then $A = \{ \emptyset \}$, $B = \{1 \}$ and $C = \{2, 3\}$. Hence, in this case the admissible subsets of $\{1,2,3\}$ are $\{1\}$, $\{2\}$, $\{3\}$ and $\{2,3\}$.  
\end{example}


The isomorphisms $N^{1}(\pe) \cong N^{1}(X) \oplus {\bf{Z}}$ and $N^{1}(X) \cong {\bf{Z}}^{d-n}$, described in \S \ref{subsection.toric} and \S \ref{subsection.vanishing.orders}, induce an isomorphism between $N^{1}(\pe)_{{\bf{R}}}$ and ${\bf{R}}^{d-n+1}$, which we use to identify these spaces hereafter. 
Likewise, we identify ${\bf{R}}^{n+1} \times N^{1}(\pe)_{{\bf{R}}}$ with ${\bf{R}}^{d+2}$, with coordinates $ ( x_{1} , \ldots ,  \linebreak[0] x_{n+1} ,  \linebreak[0] w_{n+1} , \linebreak[0] \ldots ,  \linebreak[0] w_{d},  \linebreak[0] w ) $. Let $v_{1}^{*},\ldots, v_{n}^{*}$ be the basis of $M_{{\bf{R}}} =_{def} M \otimes {\bf{R}}$ dual to the basis $v_{1},\ldots, v_{n}$ of $N_{{\bf{R}}}=_{def} N \otimes {\bf{R}}$. Note that we have an isomorphism
\begin{align*}
\begin{matrix}
\psi \colon & {\bf{R}}^{d+2} & \longrightarrow & M_{{\bf{R}}} \times {\bf{R}}^{d-n+2} \\
&(x_{1},\ldots,x_{n},x_{n+1},w_{n+1},\dots,w_{d},w) & \longmapsto & \big{(} -\sum_{i=1}^{n}{x_{i}v_{i}^{*}} + x_{n+1} u_{1} + (w-x_{n+1}) u_{2} ,  \\ 
& &  & x_{n+1}, w_{n+1}, \ldots, w_{d}, w \big{)}
\end{matrix}
\end{align*}

For each $j \in \{1,\ldots, d \}$ we define the linear function: 
\begin{align*}
\begin{matrix}
\gamma_{\mathcal{E},j} \colon& M_{{\bf{R}}} \times {\bf{R}}^{d-n+2} & \longrightarrow & {\bf{R}} \\
&(u,x_{n+1}, w_{n+1},\ldots,w_{d},w) & \longmapsto & \frac{\langle u,v_{j} \rangle -a_{j}w - w_{j}}{b_{j}-a_{j}}
\end{matrix}
\end{align*}
for any $u \in M_{{\bf{R}}}$, and any $x_{n+1},w_{n+1},\ldots,w_{d},w \in {\bf{R}}$, and where $w_{j}=0$ for each $j \leq n $. We will denote this function simply by $\gamma_{j}$, when no confusion is likely to arise.
Finally, for each admissible set $J \subseteq \{1,\ldots,d\}$, we define the linear function $I_{J} \colon {\bf{R}}^{d+2} \rightarrow {\bf{R}}$ by declaring its value at $P = (x_{1},\ldots,x_{n+1},w_{n+1},\ldots,w_{d},w) \in {\bf{R}}^{d+2}$ to be:
\begin{equation*} 
I_{J}(P)=
	\begin{cases} 
		\ \gamma_{j} \circ \psi (P),& \text{ if $J = \{ j \} \subseteq A$},\\
		\ \gamma_{j} \circ \psi (P) - x_{n+1},& \text{ if $J = \{ j \} \subseteq B$}, \\
		\ \sum_{j \in J}{\gamma_{j} \circ \psi (P) }-w+x_{n+1},& \text{ if $J \subseteq C$}.
	\end{cases}
\end{equation*}
For notational convenience, we define for each admissible set $J \subseteq \{1,\ldots, d \}$ the linear function $I'_{J} \colon M_{{\bf{R}}} \times {\bf{R}}^{d-n+2} \rightarrow {\bf{R}}$ to be $I_{J}' = I_{J} \circ \psi^{-1}$.


\subsection{The global Okounkov body of $\pe$} \label{subsection.theorem}


\begin{theorem} \label{theorem}
Let $\mathcal{E}$ be a toric vector bundle of rank two on the smooth projective toric variety $X$. 
The global Okounkov body $\Delta ({\bf{P}}(\mathcal{E}))$ of ${\bf{P}}(\mathcal{E})$ is the rational polyhedral cone in
${\bf{R}}^{n+1} \times N^{1} ({\bf{P}}(\mathcal{E}))_{{\bf{R}}} \simeq {\bf{R}}^{d+2}$ given by
	\begin{equation*}
		\begin{split}
			{\bf{\Delta}} = \Big{\{}(&x_{1},\ldots,x_{n+1},w_{n+1},\ldots,w_{d},w) \in {{\bf{R}}}^{d+2} \ | \ w \geq x_{n+1} \geq 0 
												\textnormal{ and } \\
							         &	I_{J}(x_{1},\ldots,x_{n+1},w_{n+1},\ldots,w_{d},w) \leq 0 \textnormal{ for all admissible } J \subseteq 		
												\{1,\ldots,d\} \Big{\}}.
		\end{split}
	\end{equation*}
\end{theorem}
\begin{proof}
%
%
From the characterization of the global Okounkov body in terms of its fibers over big classes in $N^{1} ({\bf{P}}(\mathcal{E}))_{{\bf{Q}}}$, it suffices to show the following stronger assertion: For every class $\mathcal{L} \in N^{1} ({\bf{P}}(\mathcal{E}))_{{\bf{Q}}}$, the fiber ${\bf{\Delta}}_{\mathcal{L}}$ of ${\bf{\Delta}}$ over $\mathcal{L}$ is equal to $\Delta(\mathcal{L}) \times \{ \mathcal{L} \}$.

To prove the assertion, we consider a class $\mathcal{L} \in N^{1} ({\bf{P}}(\mathcal{E}))_{{\bf{Q}}}$. Note that ${\bf{\Delta}}_{\mathcal{L}^{m}} = m {\bf{\Delta}}_{\mathcal{L}}$ and $\Delta(\mathcal{L}^{m}) \times \{ \mathcal{L}^{m} \} = m (\Delta(\mathcal{L}) \times \{ \mathcal{L} \})$, for each $m \in {\bf{Z}}^{+}$. Hence, we can assume that $\mathcal{L} \in N^{1} ({\bf{P}}(\mathcal{E}))$. Let $w_{n+1},\ldots,w_{d},w \in {\bf{Z}}$ be the unique integers such that 
\[
\mathcal{L} = \Ope(w)\otimes\pi^{*}\Ox \Big(\sum_{i=n+1}^{d}{w_{i}D_{i}}\Big).
\]
For notational convenience, we set $w_{i}=0$ for each $i \in \{ 1,\ldots,n \}$.
We first show that $\Delta(\mathcal{L}) \times \{ \mathcal{L} \} \subseteq {\bf{\Delta}}_{\mathcal{L}}$. For this, it is enough to see that the set
	\begin{equation*}
\Delta(\mathcal{L}) \times \{ \mathcal{L} \}	=		\overline{ \operatorname{Conv} \Big{(} \bigcup_{m \in {\bf{Z}}^{+}} \frac{1}{m} \nu \big{(} \mathcal{L}^{m} \big{)} 
				\Big{)} } \times (w_{n+1},\ldots,w_{d},w)
	\end{equation*}
is contained in ${\bf{\Delta}}$. Since ${\bf{\Delta}}$ is closed and convex, it suffices show that
\[
\big{(}\frac{1}{m}\nu(\mathcal{L}^{m})\big{)} \times (w_{n+1},\ldots,w_{d},w) \subseteq {\bf{\Delta}},
\] 
for each $m \in {\bf{Z}}^{+}$.
Furthermore, since ${\bf{\Delta}}$ is a cone, it is enough to prove that this inclusion holds when $m=1$.
With this in mind, we consider 
\[P=(x_{1},\ldots,x_{n+1},w_{n+1},\ldots,w_{d},w) \in \nu(\mathcal{L}) \times (w_{n+1},\ldots,w_{d},w).\]
Note that the existence of $P$ implies that $w \geq 0$. Let $Q = \psi(P)$, i.e.
\[Q=( - \sum_{i=1}^{n}{x_{i}v_{i}^{*}} + x_{n+1}u_{1} + (w-x_{n+1}) u_{2} , x_{n+1}, w_{n+1}, \ldots, w_{d}, w ) \in M_{{{\bf{R}}}} \times {{\bf{R}}}^{d-n+2}.\] 
By replacing $\mathcal{L}$ with a suitable tensor power, we can assume that $\gamma_{j}(Q) \in {\bf{Z}}$ for all $j \in \{1,\ldots,d \}$.
By the projection formula we have $\pi_{*}\mathcal{L}=(Sym^{w}\mathcal{E}) \otimes \Ox (\sum_{i=1}^{d}{w_{i}D_{i}})$.
If we denote the fiber of $\Ox(\sum^{d}_{i=1}{w_{i}D_{i}})$ over the unit of the torus by $G$, then by Example \ref{tensor.product} we get
\[
(\pi_{*}\mathcal{L})^{\rho_{j}}(i)=(Sym^{w}\mathcal{E})^{\rho_{j}}(i-w_{j}) \otimes G \text{,}
\] for all $j \in \{ 1,\ldots,d \}$ and all $i \in {\bf{Z}}$. Since $(x_{1},\ldots,x_{n+1}) \in \nu(\mathcal{L})$, Proposition \ref{Proposition:OnTypesOfSectionsThatAreEnough} implies that there exist $u \in M$ and a nonzero section $s \in H^{0}(\pe,\mathcal{L})_{u}=H^{0}(X,\pi_{*}\mathcal{L})_{u}$  
such that $\nu_{Y_{\bullet}}(s)=(x_{1},\ldots,x_{n+1})$, and such that
\[
s(t_{0})\in V=_{def} Sym_{E}^{w}(E_{1}^{\alpha_{0}},L_{1}^{\alpha_{1}},\ldots,L_{p}^{\alpha_{p}},L^{\alpha}) \otimes G,
\] 
for some one-dimensional subspace $L$ of $E$ different from $E_{1},L_{1},\ldots,L_{p}$ and for some $\alpha_{0},$ $ \ldots,$ $\alpha_{p},$ $\alpha \in {\bf{Z}}_{\geq 0}$ with $\sum^{p}_{i=0}{\alpha_{i}}+\alpha=w$.
Note now that
\begin{align*}
0 \neq V &\subseteq \operatorname{Im} \Big( H^{0} \big( X,\pi_{*}\mathcal{L} \big)_{u} \longhookrightarrow (Sym^{w}E) \otimes G          \Big) = \bigcap^{d}_{i=1}(\pi_{*}\mathcal{L})^{\rho_{i}}(\langle u,v_{i} \rangle) \subseteq (\pi_{*}\mathcal{L})^{\rho_{j}}(\langle u,v_{j} \rangle),
\end{align*}
for each $j \in \{ 1,\ldots,d \}$. 
By Lemma \ref{Lemma:OrdersOfVanishing}, we have that $\alpha_{0}=x_{n+1}$ and $x_{i} = \langle -u+\alpha_{0}u_{1}+(w-\alpha_{0})u_{2}, v_{j} \rangle$ for each $i \in \{ 1,\ldots,n \}$. 
In particular, we see that $Q=\left( u , x_{n+1}, w_{n+1},\ldots, w_{d}, w \right)$. 
From $\alpha_{0}=x_{n+1}$ we get
\begin{equation}\label{topandbottom}
w \geq x_{n+1} \geq 0. 
\end{equation}
Hence, we are reduced to proving that $I_{J}(P)\leq 0$ for each admissible set $J \subseteq \{1,\ldots,d\}$, or equivalently, to proving that $I_{J}'(Q)\leq 0$ for every such $J$.\\ 
Let us consider an admissible set $J \subseteq \{1,\ldots,d\}$. Then either $J=\{j\}$ for some $j \in A$, $J=\{j\}$ for some $j \in B$, or $J=\{j_{1},\ldots,j_{l}\}$ for some $j_{1},\ldots,j_{l} \in C$ such that there exist distinct indices $i_{1},\ldots,i_{l} \in \{1,\ldots,p\}$, with $j_{h} \in C_{i_{h}}$ for each $h \in \{ 1,\ldots,l \}$.\\
In the first case $0 \neq V \subseteq (\pi_{*}\mathcal{L})^{\rho_{j}}(\langle u,v_{j} \rangle)$ gives $\langle u,v_{j} \rangle \leq a_{j}w+w_{j}$, and then
\begin{equation}\label{case1}
I_{J}'(Q)\leq0.
\end{equation}
In the second case $0 \neq V \subseteq (\pi_{*}\mathcal{L})^{\rho_{j}}(\langle u,v_{j} \rangle)= Sym_{E}^{w}(E_{1}^{max\{0,\gamma_{j}(Q)\}}) \otimes G$, and this implies
	\begin{multline*}	
		0 \neq V = V \cap (Sym_{E}^{w}(E_{1}^{max\{0,\gamma_{j}(Q)\}}) \otimes G) \subseteq (Sym_{E}^{w}(E_{1}^{max\{0,\gamma_{j}(Q)\}}) \otimes G) \cap \\  (Sym_{E}^{w}(L_{1}^{\alpha_{1}}) \otimes G) \cap \cdots \cap (Sym_{E}^{w}(L_{p}^{\alpha_{p}}) \otimes G) \cap (Sym_{E}^{w}(L^{\alpha}) \otimes G),
	\end{multline*}
and then from Lemma \ref{Lemma:SubspacesOfSymE} we get $\gamma_{j}(Q) \leq max \{ 0,\gamma_{j}(Q)\} \leq w- \sum^{p}_{i=1}{\alpha_{i}}-\alpha = x_{n+1}$, so
\begin{equation}\label{case2}
I_{J}'(Q)\leq0.
\end{equation}
In the third case, we have $\gamma_{j_{h}}(Q) \leq \alpha_{j_{h}}$ for each $h \in \{ 1,\ldots,l \}$. Otherwise we would have $0 \neq V \subseteq (\pi_{*}\mathcal{L})^{\rho_{j_{h}}}(\langle u,v_{j_{h}} \rangle)= Sym_{E}^{w}( L_{i_{h}}^{\gamma_{j_{h}}(Q)} ) \otimes G$, and from Lemma \ref{Lemma:SubspacesOfSymE} we would get
\begin{multline*}
		V = V \cap (Sym_{E}^{w}( L_{i_{h}}^{\gamma_{j_{h}}(Q)} ) \otimes G) = (Sym_{E}^{w}(E_{1}^{\alpha_{0}}) \otimes G) \cap (Sym_{E}^{w}(L_{1}^{\alpha_{1}}) \otimes G) \cap \cdots \cap \\ (Sym_{E}^{w}(L_{i_{h}-1}^{\alpha_{i_{h}-1}}) \otimes G) \cap (Sym_{E}^{w}(L_{i_{h}}^{\gamma_{j_{h}(Q)}}) \otimes G) \cap (Sym_{E}^{w}(L_{i_{h}+1}^{\alpha_{i_{h}+1}}) \otimes G) \cap 
						\cdots \cap \\ (Sym_{E}^{w}(L_{p}^{\alpha_{p}}) \otimes G) \cap (Sym_{E}^{w}(L^{\alpha}) \otimes G)	= 0,
\end{multline*}
which is a contradiction. By adding these inequalities over $j_{h} \in J$, we get
\[
\sum_{j\in J}{\gamma_{j}(Q)}\leq \sum^{l}_{h=1}{\alpha_{j_{h}}}\leq \sum^{p}_{i=1}{\alpha_{i}}=w-\alpha_{0}-\alpha\leq w-x_{n+1},
\]
and therefore 
\begin{equation}\label{case3}
I_{J}'(Q)\leq0.
\end{equation}
From (\ref{topandbottom}), (\ref{case1}), (\ref{case2}) and (\ref{case3}), it follows that $P\in {\bf{\Delta}}$.

%
%
As we have completed the proof of $\Delta(\mathcal{L}) \times \{ \mathcal{L} \} \subseteq {\bf{\Delta}}_{\mathcal{L}}$, we now prove that ${\bf{\Delta}}_{\mathcal{L}} \subseteq \Delta(\mathcal{L}) \times \{ \mathcal{L} \}$. 
For this, we note that
\[
{\bf{\Delta}}_{\mathcal{L}}=\overline{{\bf{\Delta}}_{\mathcal{L}} \cap \left( {{\bf{Q}}}^{n+1} \times N^{1}(\pe)_{{\bf{Q}}} \right)},
\] 
since ${\bf{\Delta}}_{\mathcal{L}}$ is defined by rational linear inequalities. Thus, it suffices to show that
\[
{\bf{\Delta}}_{\mathcal{L}} \cap \left( {{\bf{Q}}}^{n+1} \times N^{1}(\pe)_{{\bf{Q}}} \right) \subseteq \Delta(\mathcal{L}) \times \{ \mathcal{L} \}. 
\]
To prove this, let us consider
\[ 
P=(x_{1},\ldots,x_{n+1},w_{n+1},\ldots,w_{d},w) \in {\bf{\Delta}}_{\mathcal{L}} \cap \left( {\bf{Q}}^{n+1} \times N^{1}(\pe)_{{\bf{Q}}} \right),
\]
and define $Q$ to be $\psi(P)$, i.e. 
\[
Q=(-\sum_{i=1}^{n}x_{i}v_{i}^{*}+x_{n+1}u_{1}+(w-x_{n+1})u_{2},x_{n+1},w_{n+1},\ldots,w_{d},w) \in M_{{\bf{Q}}} \times {\bf{Q}} \times N^{1}(\pe)_{{\bf{Q}}}.
\]
By replacing $\mathcal{L}$ with a suitable tensor power, we can assume that $Q \in M \times {\bf{Z}} \times N^{1}(\pe)$ and $\gamma_{j}(Q) \in {\bf{Z}}$ for each $j \in \{ 1,\ldots, d \}$.
We have $w \geq x_{n+1} \geq 0$, and $I_{J}'(Q) = I_{J}(P) \leq 0$ for each admissible set $J \subseteq \{1,\ldots,d\}$. Let us define 
\[
u = -  \sum_{i=1}^{n}x_{i}v_{i}^{*}+  x_{n+1}u_{1} +  (w-x_{n+1})u_{2} \in M,
\]
and note that $Q = ( u ,x_{n+1},w_{n+1},\ldots,w_{d},w)$.
We will show the existence of a nonzero section $s \in H^{0}(\pe,\mathcal{L})_{u} = H^{0} (X,\pi_{*}\mathcal{L})_{u}$,
satisfying $\nu_{\bullet}(s) = ( x_{1} , \ldots , x_{n+1})$, which will give us
\[
P\in  \nu (\mathcal{L}) \times (w_{n+1},\ldots,w_{d},w) \subseteq \Delta(\mathcal{L}) \times \{ \mathcal{L} \}.
\] 
By the projection formula we have $\pi_{*}\mathcal{L}=(Sym^{w}\mathcal{E}) \otimes \Ox (\sum_{i=1}^{d}{w_{i}D_{i}})$.
If we denote the fiber of $\Ox(\sum^{d}_{i=1}{w_{i}D_{i}})$ over the unit of the torus by $G$, then by Example \ref{tensor.product} we get
\[
(\pi_{*}\mathcal{L})^{\rho_{j}}(i) = (Sym^{w}\mathcal{E})^{\rho_{j}}(i-w_{j}) \otimes G \text{,}
\] for all $j \in \{ 1,\ldots,d \}$ and all $i \in {\bf{Z}}$. 
For each $i \in \{ 1,\ldots,p \}$, define $\alpha_{i}\in {\bf{Z}}_{\geq 0}$ by
\[
\alpha_{i}=_{def} max \left(\{0\} \cup \{\gamma_{j}(Q) \ | \ j\in C_{i} \} \right).
\]
We claim that $\alpha=_{def} w-x_{n+1}-\sum^{p}_{i=1}{\alpha_{i}}$ is a nonnegative integer. Indeed, this is clear if $\alpha_{i}=0$ for each $i \in \{1,\ldots , p \}$. On the other hand, if $\alpha_{i} \neq 0$ for some $i \in \{1,\ldots , p \}$, let $i_{1},\ldots,i_{l} \in \{ 1,\ldots,p \}$ be the distinct indices such that for $i \in \{ 1,\ldots,p \} $, $\alpha_{i}\neq 0$ if and only if $i \in \{i_{1}, \ldots, i_{l}\}$. For each $h \in \{ 1,\ldots,l \}$, let us choose $j_{h}\in \{ 1,\ldots,d \}$ such that $j_{h}\in C_{i_{h}}$ and $\alpha_{i_{h}}=\gamma_{j_{h}}(Q)$. Then the set $J=\{ j_{1}\ldots,j_{l} \}$ is admissible, and we have $I_{J}'(Q)\leq 0$. Hence 
\[
\sum^{p}_{i=1}{\alpha_{i}} = \sum^{l}_{h=1}{\alpha_{i_{h}}} = \sum_{j\in J}{\gamma_{j}(Q)} \leq w-x_{n+1}.
\] 
In either case, it follows that $\alpha$ is a nonnegative integer. Let $L$ be a one-dimensional subspace of $E$ different from $E_{1},L_{1},\ldots,L_{p}$. From Lemma \ref{Lemma:SubspacesOfSymE}, we see that
\[
V=_{def} Sym_{E}^{w}(E_{1}^{x_{n+1}},L_{1}^{\alpha_{1}},\ldots,L_{p}^{\alpha_{p}},L^{\alpha}) \otimes G
\]
is a one-dimensional subspace of $(Sym^{w}E) \otimes G$.
We now prove that $V \subseteq (\pi_{*}\mathcal{L})^{\rho_{j}}(\langle u,v_{j} \rangle)$ for each $j \in \{ 1,\ldots,d \}$, considering separately the cases $j\in A$, $j\in B$ and $j\in C$.
If $j\in A$, then $J=\{j\}$ is admissible, and $I_{J}'(Q) \leq 0$. This gives $\langle u,v_{j}\rangle \leq a_{j}w+w_{j}$, and therefore
\[
V \subseteq (Sym^{w}E) \otimes G = (\pi_{*}\mathcal{L})^{\rho_{j}}(\langle u,v_{j} \rangle).
\] 
If $j\in B$, then $J=\{j\}$ is admissible, and $I_{J}'(Q) \leq 0$. This gives 
\[
\langle u,v_{j}\rangle \leq a_{j}w+w_{j}+(b_{j}-a_{j})x_{n+1}, 
\]
and therefore
\[
V \subseteq Sym_{E}^{w}(E_{1}^{x_{n+1}}) \otimes G = (\pi_{*}\mathcal{L})^{\rho_{j}}(a_{j}w+w_{j}+(b_{j}-a_{j})x_{n+1}) \subseteq (\pi_{*}\mathcal{L})^{\rho_{j}}(\langle u,v_{j} \rangle).
\]
If $j\in C$, then there exists $i\in \{ 1,\ldots,p \}$ such that $j\in C_{i}$, and
\[
V \subseteq Sym_{E}^{w}(L_{i}^{\alpha_{i}}) \otimes G \subseteq Sym_{E}^{w}(L_{i}^{max\{ 0,\gamma_{j}(Q)\}}) \otimes G = (\pi_{*}\mathcal{L})^{\rho_{j}}(\langle u,v_{j} \rangle).
\]
Therefore
\[
V \subseteq \bigcap^{d}_{j=1}(\pi_{*}\mathcal{L})^{\rho_{j}}(\langle u,v_{j} \rangle)=\operatorname{Im} \Big( H^{0} \big( X,\pi_{*}\mathcal{L}) \big)_{u} \longhookrightarrow (Sym^{w}E) \otimes G\Big).
\]
We can now choose a nonzero section $s \in H^{0}(\pe,\mathcal{L})_{u}=H^{0}(X,\pi_{*}\mathcal{L})_{u}$ such that $s(t_{0}) \in V$. By Lemma \ref{Lemma:OrdersOfVanishing}, the section $s$ satisfies 
\begin{align*}
\nu_{Y_{\bullet}}(s) &= \big( \langle -u+  x_{n+1} u_{1}+(w-  x_{n+1} )u_{2}, v_{1} \rangle ,\ldots, \langle-u+  x_{n+1} u_{1}+(  w -   x_{n+1} ) u_{2}, v_{n} \rangle ,  x_{n+1} \big)  \\ &= (x_{1},\ldots,x_{n+1}).
\end{align*}
Thus $P=\nu_{Y_{\bullet}}(s) \times (w_{n+1},\ldots,w_{d},w) \in \Delta(\mathcal{L}) \times \{ \mathcal{L} \}$.
It follows that ${\bf{\Delta}}_{\mathcal{L}} \subseteq \Delta(\mathcal{L}) \times \{ \mathcal{L} \}$, and 
this completes the proof.
\end{proof}




\begin{remark} \label{explicit}Explicitly, the inequalities defining the global Okounkov body $\Delta ({\bf{P}}(\mathcal{E}))$ of ${\bf{P}}(\mathcal{E})$ given in Theorem \ref{theorem} are $w \geq x_{n+1} \geq 0$ and:
\begin{equation*}
\begin{split}
& \sum_{i=1}^{n} \langle v_{i}^{*}, v_{j} \rangle x_{i} + \langle u_{2} - u_{1} , v_{j} \rangle x_{n+1} +( a_{j} - \langle u_{2} , v_{j} \rangle ) w + w _{j} \geq 0, \text{ for each } j \in A, \\
& \sum_{i=1}^{n} \langle v_{i}^{*}, v_{j} \rangle x_{i} + (\langle u_{2} - u_{1} , v_{j} \rangle + b_{j} - a_{j} ) x_{n+1} + ( a_{j} -\langle u_{2} , v_{j} \rangle ) w + w _{j} \geq 0, \text{ for each } j \in B, \\
& \sum_{j \in J} \frac{1}{b_{j}-a_{j}} \left[ \sum_{i=1}^{n} \langle v_{i}^{*}, v_{j} \rangle x_{i} + \langle u_{2} - u_{1} , v_{j} \rangle  x_{n+1} + ( a_{j} - \langle u_{2} , v_{j} \rangle ) w + w _{j}   \right] + w - x_{n+1} \geq 0, \\
& \text{ for each admissible set } J \subseteq C. \\
\end{split}
\end{equation*}
\end{remark}



\begin{application} The next proposition shows that the projectivization of a rank two toric vector bundle over a simplicial projective toric variety is a Mori dream space. This result has been proven in \cite{Knop} and \cite{Suess} in the more general setting of $T$-varieties of complexity one, i.e. normal varieties with an algebraic action of a torus $T$ such that the lowest codimension of a $T$-orbit is equal to one. Hu and Keel introduced Mori dream spaces in \cite{HK} as a class of varieties with interesting features from the point of view of Mori theory, for example, on these varieties the Mori program can be carried out for any pseudoeffective divisor. A projective ${\bf{Q}}$-factorial variety $Z$ with $\operatorname{Pic}(Z)_{{\bf{Q}}}=N^{1}(Z)_{{\bf{Q}}}$ is a Mori dream space if it has a finitely generated Cox ring (see Proposition 2.9 in \cite{HK}). This use of Okounkov bodies could provide new insights into investigating the finite generation of Cox rings of different varieties, for instance of projectivizations of higher rank toric vector bundles (see Question 7.2 in \cite{PTVB}).

\begin{proposition} \label{theorem2}
Any Cox ring of the projectivization $\pe$ of a rank two toric vector bundle $\mathcal{E}$ over the projective simplicial toric variety $X$ is finitely generated and $\pe$ is a Mori dream space.
\end{proposition}
\begin{proof}
Any simplicial toric variety is ${\bf{Q}}$-factorial, and a projective bundle over a ${\bf{Q}}$-factorial variety is again ${\bf{Q}}$-factorial, hence the we are reduced to prove the finite generation of any Cox ring of $\pe$ in the sense of Hu and Keel. We consider a toric resolution of singularities $f \colon X' \rightarrow X$, i.e. $X'$ is a smooth toric variety and $f$ is a proper birational toric morphism. Given a toric vector bundle $\mathcal{E}$ on $X$, the induced map $f' \colon {\bf{P}}(f^{*}{\mathcal{E}}) \rightarrow \pe$ is also proper and birational. In this case $f^{*}{\mathcal{E}}$ is a toric vector bundle on $X'$ and the finite generation of a Cox ring of ${\bf{P}}(f^{*}{\mathcal{E}})$ implies the finite generation of any Cox ring of $\pe$. Therefore we can assume that the toric variety $X$ is smooth and projective. Let us prove that the semigroup defined by
\[
\begin{split}
S = \bigl\{( & x_{1},\ldots, x_{n+1},w_{n+1},\ldots,w_{d},w) \in {\bf{R}}^{d+2} = {\bf{R}}^{n+1} \times N^{1}(\pe)_{{\bf{R}}} \ | \ \textnormal{ There exists } \\
& s \in H^{0} \bigl( \pe, \Ope(w) \otimes \pi^{*}\Ox(\sum^{d}_{i=n+1}w_{i}D_{i}) \bigr) \textnormal{ such that } \nu_{Y_{\bullet}}(s) = (x_{1},\ldots,x_{n+1}) \bigr\}
\end{split}
\]
is finitely generated. Since $S \subseteq {\bf{Z}}^{d+2}$, it is enough to prove that the semigroup $S \cap (c \cdot {\bf{Z}}^{d+2})$ is finitely generated, where $c = \operatorname{lcm} \{ b_{j}-a_{j} \, | \, j=1,2,\ldots,d \}$. And for this it suffices to prove that $S \cap (c \cdot {\bf{Z}}^{d+2}) = \Delta(\pe) \cap (c \cdot {\bf{Z}}^{d+2})$, since $\Delta(\pe)$ is a rational polyhedral cone. From the definition of $\Delta(\pe)$, we have that $S \cap (c \cdot {\bf{Z}}^{d+2}) \subseteq \Delta(\pe) \cap (c \cdot {\bf{Z}}^{d+2})$. Let $( x_{1},\ldots, x_{n+1},w_{n+1},\ldots,w_{d},w) \in \Delta(\pe) \cap (c \cdot {\bf{Z}}^{d+2})$. As in the second part of the proof of Theorem \ref{theorem}, it follows that there exists a nonzero section $s \in H^{0} \bigl( \pe, \Ope(w) \otimes \pi^{*}\Ox(\sum^{d}_{i=n+1}w_{i}D_{i}) \bigr)$ that satisfies $\nu_{Y_{\bullet}}(s) = (x_{1},\ldots,x_{n+1})$. Therefore the semigroup $S$ is finitely generated as we claimed. Now we prove that the Cox ring of $\pe$ associated to the line bundles $\pi^{*}\Ox(D_{n+1}),\ldots,\pi^{*}\Ox(D_{d})$ and $\Ope(1))$ on $\pe$ 
is finitely generated. This Cox ring is equal to 
\[
R = \bigoplus_{(m_{n+1},\ldots,m_{d},m) \in {\bf{Z}}^{d-n+1}} R_{(m_{n+1},\ldots,m_{d},m)},
\]
where for each $(m_{n+1},\ldots,m_{d},m) \in {\bf{Z}}^{d-n+1}$,
\[
R_{(m_{n+1},\ldots,m_{d},m)} \ =_{def}  H^{0} (X, (Sym^{m} \mathcal{E}) \otimes  \Ox(m_{n+1}D_{n+1}) \otimes \cdots \otimes \Ox(m_{d}D_{d})).
\]
Let $\{g_{1},g_{2},\ldots,g_{l}\}$ be generators of $S$. 
For each $j \in \{1,2,\ldots,l \}$, there exist $m_{n+1}^{(j)}, \linebreak[0] \ldots , \linebreak[0] m_{d}^{(j)}, \linebreak[0] m^{(j)} \in {\bf{Z}}$ and a nonzero section $s_{j} \in H^{0} \bigl( \pe, \Ope(m^{(j)}) \otimes \pi^{*}\Ox(\sum^{d}_{i=n+1}m_{i}^{(j)}D_{i}) \bigr)$ such that $g_{l} = ( \nu_{Y_{\bullet}}(s_{j}), m_{n+1}^{(j)},\ldots, m_{d}^{(j)}, m^{(j)} )$.
For each $(m_{n+1}, \linebreak[0] \ldots , \linebreak[0] m_{d} , \linebreak[0] m) \in {\bf{Z}}^{d-n+1}$, since $g_{1},g_{2},\ldots,g_{l}$ generate $S$, it follows that
\[
\nu_{Y_{\bullet}}\bigl( ( k[s_{1},s_{2},\ldots,s_{l}] \cap   R_{(m_{n+1},\ldots,m_{d},m)} )  \smallsetminus \{ 0 \} \bigr) = \nu_{Y_{\bullet}}\bigl(R_{(m_{n+1},\ldots,m_{d},m)} \smallsetminus \{ 0 \} \bigr).
\]
By Remark \ref{Lemma:LazMus} the vector spaces $k[s_{1},s_{2},\ldots,s_{l}] \cap   R_{(m_{n+1},\ldots,m_{d},m)}$ and $R_{(m_{n+1},\ldots,m_{d},m)}$ have the same dimension, and thus they are equal. Therefore $R = k[s_{1},s_{2},\ldots,s_{l}]$ and this completes the proof.
\end{proof}

\end{application}



\section{Examples} \label{section.examples}
The explicit description of Okounkov bodies in concrete examples can be rather difficult. Our result allows us to explicitly compute the Okounkov bodies of all line bundles on projectivizations of rank two toric vector bundles over smooth projective toric varieties, with respect to the flag from \S \ref{flag}, by substituting combinatorial data into the inequalities given in Remark \ref{explicit}. In this section we present some examples to illustrate our main result.

\begin{example}
We consider $T_{{\bf{P}}^{2}}$, the tangent bundle of the projective plane (see Example \ref{tangent.bundle.definition}). From Remark \ref{explicit} (see Example \ref{tangent.bundle.definition.2}), we get inequalities for the Okounkov body of each line bundle on ${\bf{P}}(T_{{\bf{P}}^{2}})$. For instance, by setting $w=1$ and $w_{j}=0$ for each $j$, we deduce that the Okounkov body $\Delta( \mathcal{O}_{{\bf{P}}(T_{{\bf{P}}^{2}})}(1))$ is defined inside ${\bf{R}}^{3}$ by the inequalities: 
\begin{align*}
\begin{matrix}
 1 \geq x_{3} & & & x_{3} \geq 0 & & & x_{1} \geq 0	\\
 x_{2} \geq 0 & & & 2 \geq x_{1}+x_{2}+x_{3} & & & 1 \geq x_{1} \\
\end{matrix}
\end{align*}
In particular, we see that the volume of $\mathcal{O}_{{\bf{P}}(T_{{\bf{P}}^{2}})}(1)$ is $\textnormal{vol}_{{\bf{R}}^3}\big( \Delta(\mathcal{O}_{{\bf{P}}(T_{{\bf{P}}^{2}})}(1)) \big) \cdot 3! = 6$.


\begin{center}
\includegraphics[width=3.7in]{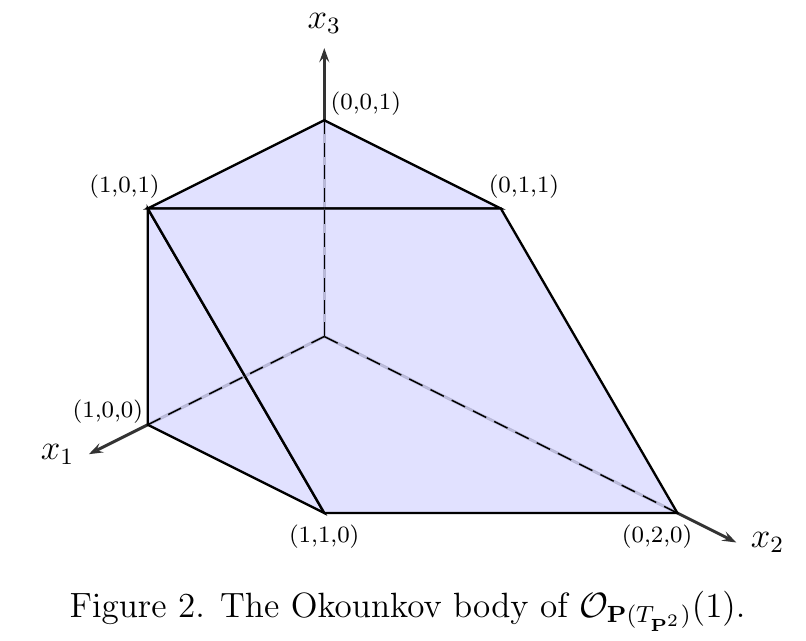}
\end{center}
\end{example}


In the next example we see that our description gives the expected answer for line bundles that are pulled back from the base.

\begin{example}
The inequalities for the Okounkov body of a line bundle on $\pe$ of the form $\pi^{*}\mathcal{O}_{X}(\sum_{i=n+1}^{d}m_{i}D_{i})$ are $x_{n+1}=0$ and
\[
\sum_{i=1}^{n}\langle v^{*}_{i}, v_{j} \rangle x_{i} + m_{j} \geq 0,
\]
for each $j \in \{1,\ldots,d\} $.
Furthermore, from the description of the Okounkov body of a toric line bundle on $X$ given in \cite[Proposition 6.1]{CBALS}, we see that 
\[
\Delta_{Y_{\bullet}}\bigl(\pi^{*}\mathcal{O}_{X} ( \sum_{i=n+1}^{d}m_{i}D_{i} ) \bigr) =
P_{\sum_{i=n+1}^{d}m_{i}D_{i}} \times \{ 0 \}
= \Delta_{X_{\bullet}}\bigl( \mathcal{O}_{X} (\sum_{i=n+1}^{d}m_{i}D_{i} ) \bigr) \times \{ 0 \}.
\]
\end{example}


In the next example, we see that our description gives the expected answer when the toric vector bundle equivariantly splits.

\begin{example}
When $\mathcal{E}$ equivariantly splits as the sum of two toric line bundles $\mathcal{L}_{1}$ and $\mathcal{L}_{2}$, the variety $\pe$ is a toric variety. The subvarieties in our flag in $\pe$ from \S \ref{flag} are also invariant with respect to the torus $T'$ of $\pe$. Hence, in this case we have two descriptions of the Okounkov bodies of line bundles on $\pe$ with respect to this flag of invariant subvarieties, namely, the one given by our theorem, and the one given by Lazarsfeld and Musta\c{t}\u{a} in \cite{CBALS} in the case of toric varieties. It is good to see that these two descriptions agree, as expected. 

Let $h_{1}$ and $h_{2}$ be the piecewise linear functions associated to $\mathcal{L}_{1}$ and $\mathcal{L}_{2}$ (see 3.4 in \cite{Fulton}). 
In particular, $\mathcal{L}_{1}$ and $\mathcal{L}_{2}$ correspond to the $T$-invariant divisors $-\sum_{j=1}^{d}h_{1}(v_{j})D_{j}$ and $-\sum_{j=1}^{d}h_{2}(v_{j})D_{j}$. 
Let $\Phi \colon N_{{\bf{R}}} \rightarrow N_{{\bf{R}}} \times {\bf{R}}$ be the piecewise linear map defined by $\Phi(v) = (v,h_{1}(v)-h_{2}(v))$. For each cone $\sigma \in \Delta$, let $\sigma^{+}$ and $\sigma^{-}$ be the cones in $N_{{\bf{R}}} \times {\bf{R}}$ spanned by $\Phi(\sigma)$ and $(0,1)$, and by $\Phi(\sigma)$ and $(0,-1)$, respectively. Let $\widetilde{\Delta}$ be the fan in $N_{{\bf{R}}} \times {\bf{R}}$ consisting of the faces of $\sigma^{+}$ and $\sigma^{-}$ for all $\sigma \in \Delta$. The toric variety associated to the fan $\widetilde{\Delta}$ is isomorphic to $\pe$ (see \S 7 in \cite{Oda2}). The rays of $\widetilde{\Delta}$ are $\rho^{+}={\bf{R}}_{\geq 0} \cdot (0,1)$, $\rho^{-}={\bf{R}}_{\geq 0} \cdot (0,-1)$ and $\widetilde{\rho_{j}} = \Phi(\rho_{j})$ for each $j \in \{ 1, \ldots, d \}$. Let us denote the corresponding $T'$-invariant divisors by $D^{+}$, $D^{-}$ and $\widetilde{D_{j}}$ for each $j \in \{ 1, \ldots, d \}$, respectively. The flag in $\pe$ given by our construction from \S \ref{flag} is 
\[
Y_{\bullet} \colon \pe \supseteq \widetilde{D_{1}} \supseteq \widetilde{D_{1}} \cap \widetilde{D_{2}} \supseteq \cdots \supseteq \widetilde{D_{1}} \cap \widetilde{D_{2}} \cap \cdots \cap \widetilde{D_{n}} \supseteq \widetilde{D_{1}} \cap \widetilde{D_{2}} \cap \cdots \cap \widetilde{D_{n}} \cap D^{-}.
\]
Note that $\Phi(v_{1}),\ldots,\Phi(v_{n})$ and $(0,-1)$ span the maximal cone $\tau^{-} \in \widetilde{\Delta}$. Let us change the reference ordered basis in $N_{{\bf{R}}} \times {\bf{R}}$ to $\{ \Phi(v_{1}),\ldots,\Phi(v_{n}), (0,-1) \}$. In these new coordinates the rays are given by $\rho^{+}={\bf{R}}_{\geq 0} \cdot (0,-1)$, $\rho^{-}={\bf{R}}_{\geq 0} \cdot (0,1)$ and $\widetilde{\rho_{j}} = {\bf{R}}_{\geq 0} \cdot (v_{j},h_{2}(v_{j}) - h_{1}(v_{j}) + \langle u_{2} ,v_{j} \rangle - \langle u_{1} , v_{j} \rangle)$ for each $j \in \{ 1, \ldots, d \}$. Using the argument in the proof of Proposition \ref{Lemma:EquivariantSectionsAreEnough}, we see that $\Ope (1) = \Ope (D^{+}) \otimes \pi^{*}\mathcal{L}_{2} = \Ope (D^{-}) \otimes \pi^{*}\mathcal{L}_{1}$. We set $D = D^{+} + \sum_{j=1}^{d}(-h_{2}(v_{j}) - \langle u_{2} ,v_{j} \rangle ) \widetilde{D_{j}}$, and note that the $T'$-invariant divisor $D$ satisfies $\Ope (1)= \Ope (D)$ and $D|_{U_{\tau^{-}}} = 0$.

Let us consider a line bundle $\mathcal{L} = \Ope (m) \otimes \pi^{*} \mathcal{O}_{X}(\sum_{i=n+1}^{d}m_{i}D_{i})$ on $\pe$. Let us identify the dual of $N_{{\bf{R}}} \times {\bf{R}}$ with ${\bf{R}}^{n+1}$ by identifying the ordered basis $\{ \Phi(v_{1}),\ldots,\Phi(v_{n}), (0,-1) \}$ of $N_{{\bf{R}}} \times {\bf{R}}$ with the coordinates $x_{1},\ldots,x_{n+1}$ on ${\bf{R}}^{n+1}$. We set $m_{i}=0$ for each $i \in \{ 1,\ldots,n \}$. On the one hand, the description in \cite[Proposition 6.1]{CBALS} says that with this identification $\Delta_{Y_{\bullet}}(\mathcal{L})$ is the polytope $P_{mD+\sum_{i=n+1}^{d}m_{i}\widetilde{D_{i}}}$. This polytope is defined as a subset of ${\bf{R}}^{n+1}$ by the inequalities
\begin{equation*}
\begin{cases}
x_{n+1} \geq 0, \ m \geq x_{n+1}, \\
\sum_{i=1}^{n}\langle v_{i}^{*},v_{j} \rangle x_{i} + (h_{2}(v_{j}) - h_{1}(v_{j}) + \langle u_{2} ,v_{j} \rangle - \langle u_{1} , v_{j} \rangle)x_{n+1} -m h_{2}(v_{j}) \\
\qquad \quad \qquad \qquad \qquad - m \langle u_{2} , v_{j} \rangle + m_{j} \geq 0, \qquad \qquad \text{for each $j \in \{1, \ldots, d\}$}.
\end{cases}
\end{equation*}
On the other hand, the admissible subsets of $\{1,\ldots,d \}$ associated to a toric vector bundle that equivariantly splits are exactly the singletons. From Remark \ref{explicit}, our inequalities for $\Delta_{Y_{\bullet}}(\mathcal{L})$ are
\begin{equation*}
\begin{cases}
m \geq x_{n+1} \geq 0, \\
\sum_{i=1}^{n}\langle v_{i}^{*},v_{j} \rangle x_{i} + \langle u_{2} - u_{1} ,v_{j} \rangle x_{n+1} + (h_{2}(v_{j}) - h_{1}(v_{j})) x_{n+1} -m h_{2}(v_{j}) \\
\qquad \quad \qquad \qquad \qquad - m \langle u_{2} , v_{j} \rangle + m_{j} \geq 0, \qquad \qquad \text{for each $j \in \{1, \ldots, d\}$}.
\end{cases}
\end{equation*}
Therefore, the two descriptions of the Okounkov body $\Delta_{Y_{\bullet}}(\mathcal{L})$ coincide.
\end{example}





\begin{thebibliography}{99}


\bibitem{Fulton}
W. Fulton, Introduction to toric varieties, Annals of Math. Studies, Vol. 131, Princeton Univ. Press, Princeton, 1993.

\sbl

\bibitem{Fulton.Harris}
W. Fulton, J. Harris, Representation theory: a first course, Springer-Verlag, New York, 1991, Graduate Texts in Mathematics, No. 129.

\sbl

\bibitem{Hartshorne} 
R. Hartshorne, Algebraic Geometry, Springer-Verlag, New York, 1977, Graduate Texts in Mathematics, No. 52.

\sbl

\bibitem{Suess} J. Hausen, H. S\"{u}\ss, The Cox ring of an algebraic variety with torus action, preprint, arXiv:0903.4789.

\sbl

\bibitem{PTVB}
M. Hering, S. Payne, M. Musta\c{t}\u{a}, \ Positivity for toric vector bundles, \ preprint, arXiv:0805.4035.

\sbl

\bibitem{HK}
Y. Hu, S. Keel, Mori dream spaces and GIT, Michigan Math. J. 48 (2000) 331-348.

\sbl

\bibitem{KK} K. Kaveh, A. Khovanskii, Convex bodies and algebraic equations on affine varieties, preprint, arXiv:0804.4095. 

\sbl

\bibitem{Klyachko} A. Klyachko, Equivariant vector bundles on toral varieties, Math. USSR-Izv. 35 (1990), 337-375.

\sbl

\bibitem{Knop} F. Knop, \"{U}ber Hilberts vierzehntes Problem f\"{u}r Variet\"{a}ten mit Kompliziertheit eins, Math. Z. 213 (1993) 33-35.

\sbl

\bibitem{LePotier}
J. Le Potier, Lectures on vector bundles, Cambridge studies on advanced mathematics, No. 54, Cambridge university press, Cambridge, 1997.

\sbl

\bibitem{CBALS}
R. Lazarsfeld, Mircea Musta\c{t}\u{a}, Convex bodies associated to linear series, preprint, arXiv:0805.4559.

\sbl

\bibitem{PAG}
R. Lazarsfeld, Positivity in Algebraic Geometry, I \& II, Ergebnisse der Mathematik und ihrer Grenzgebiete, Vols. 48 \& 49, Springer Verlag, Berlin, 2004.

\sbl

\bibitem{Oda2}
Tadao Oda, Lectures on torus embeddings and applications, Tata Inst. Fund. Research, Bombay, No. 58, Springer-Verlag, Berlin-Heidelberg-New York, 1978.

\sbl

\bibitem{Okounkov96} A. Okounkov, Brunn-Minkowski inequality for multiplicities, Invent. Math. 125 (1996) 405-411.

\sbl

\bibitem{Okounkov03}
A. Okounkov, Why would multiplicities be log-concave?, The orbit method in geometry and physics, Progr. Math. 213 (2003) 329-347.

\sbl

\bibitem{Payne}
Sam Payne, Toric vector bundles, branched covers of fans, and the resolution property, J.~Alg. Geom. 18 (2009) 1-36.

\sbl

\bibitem{Payne2}
Sam Payne, Moduli of toric vector bundles, Compositio Math. 144 (2008) 
1199-1213.

\sbl

\end{thebibliography}
\end{document}